\input amstex
\documentstyle{amsppt}
%
%
\nopagenumbers
\accentedsymbol\tx{\tilde x}
\def\tr{\operatorname{tr}}
\def\negskp{\hskip -2pt}
\def\compos{\,\raise 1pt\hbox{$\sssize\circ$} \,}
\pagewidth{360pt}
\pageheight{606pt}
\leftheadtext{Ruslan A. Sharipov}
\rightheadtext{V-representation for normality equations \dots}
\topmatter
\title V-representation for normality equations\\
\lowercase{in geometry of generalized \uppercase{L}egendre 
transformation.}
\endtitle
\author
R.~A.~Sharipov
\endauthor
\abstract
Normality equations describe Newtonian dynamical systems admitting
normal shift of hypersurfaces. These equations were first derived
in Euclidean geo\-metry. Then very soon they were rederived in
Riemannian and in Finslerian geometry. Recently I have found that
normality equations can be derived in geometry given by classical
and/or generalized Legendre transformation. However, in this case
they appear to be written in $\bold p$-representation, i\.\,e\. in
terms of momentum covector and its components. The goal of present
paper is to transform normality equations back to
$\bold v$-representation, which is more natural for Newtonian
dynamical systems.
\endabstract
\address Rabochaya street 5, 450003, Ufa, Russia
\endaddress
\email \vtop to 20pt{\hsize=280pt\noindent
R\_\hskip 1pt Sharipov\@ic.bashedu.ru\newline
r-sharipov\@mail.ru\vss}
\endemail
\urladdr
http:/\negskp/www.geocities.com/r-sharipov
\endurladdr
\endtopmatter
\loadbold
\TagsOnRight
\document
\head
1. Newtonian dynamical systems\\
in velocity and momentum representations.
\endhead
   Traditionally first section of a paper is introduction. However
for this paper it would be rather large. The matter is that theory
of Newtonian dynamical systems admitting normal shift was discovered
and was initially developed by me and my students (see papers
\cite{1--14} and theses \cite{15} and \cite{16}). Last few years it
is developed by me alone (see papers \cite{17--26}). In order to make
easy reading my paper for people incognizant of this theory I have
broken introductory information into several parts, each with its own
title. In sections~1--5 below I give all necessary definitions along
with historical overview, and I resume in brief results of previous
paper \cite{26}.\par
    Newtonian dynamical system describes the motion of mass point with
unit mass according to Newton's second law. It is given by ordinary
differential equations 
$$
\xalignat 2
&\hskip -2em
\dot x^i=v^i,
&&\dot v^i=\Phi^i,
\tag1.1
\endxalignat
$$
where $\Phi^i=\Phi^i(x^1,\ldots,x^n,v^1,\ldots,v^n)$. In geometric
interpretation of differential equations \thetag{1.1} variables
$x^1,\,\ldots,\,x^n$ are interpreted as coordinates of moving point
$p=p(t)$ of some $n$-dimensional manifold $M$ in a local chart
of this manifold. Variables $v^1,\,\ldots,\,v^n$ are components of
velocity vector $\bold v$ tangent to $M$ at the point $p=p(t)$.
Pair $q=(p,\bold v)$ composed by a point $p\in M$ and by some tangent
vector $\bold v\in T_p(M)$ is a point of tangent bundle $TM$. In such
interpretation dynamical system \thetag{1.1} describes not only a
moving mass point, but arbitrary mechanical system with holonomic
constraints, i\.\,e\. any machine with numerous moving parts joined
by cardan joints, thumbscrews, tooth gearings and so on (see details
in Chapter~\uppercase\expandafter{\romannumeral 2} of thesis
\cite{15}). In mechanics base manifold $M$ is called {\bf configuration
space}, while tangent bundle $TM$ is called {\bf phase space}.\par
    In some cases cotangent bundle $T^*\!M$ is used as phase space
of dynamical system. Its point $q=(p,\bold p)$ is a pair of point
$p\in M$ and momentum covector $\bold p\in T^*_p(M)$. Hamiltonian
dynamical systems form an example:
$$
\xalignat 2
&\hskip -2em
\dot x^i=\frac{\partial H}{\partial p_i},
&&\dot p_i=-\frac{\partial H}{\partial x^i},
\tag1.2
\endxalignat
$$
Here $H=H(x^1,\ldots,x^n,p_1,\ldots,p_n)$ is Hamilton function.
Each Hamiltonian dynamical system \thetag{1.2} is related to some
Lagrangian dynamical system
$$
\xalignat 2
&\hskip -2em
\dot x^i=v^i,
&&\frac{d}{dt}\!\left(\frac{\partial L}
{\partial v^i}\right)=\frac{\partial L}{\partial x^i}
\tag1.3
\endxalignat
$$
with lagrange function $L=L(x^1,\ldots,x^n,v^1,\ldots,v^n)$. Dynamical
systems \thetag{1.2} and \thetag{1.3} are related to each other by means
of Legendre transformation $\lambda$ and its inverse map $\lambda^{-1}$
(see book \cite{27} for more details):
$$
\xalignat 2
&\hskip -2em
\lambda\!:\ TM\to T^*\!M,
&&\lambda^{-1}\!:\ T^*\!M\to TM.
\tag1.4
\endxalignat
$$
Nonlinear fiber-preserving maps \thetag{1.4} are given by formulas
$$
\xalignat 2
&\hskip -2em
p_i=\frac{\partial L}{\partial v^i},
&&v^i=\frac{\partial H}{\partial p_i}.
\tag1.5
\endxalignat
$$\par
   Note that if Legendre transformation \thetag{1.5} is diffeomorphic
(at least locally), then \thetag{1.3} is a special case of Newtonian
dynamical system \thetag{1.1}. Indeed, functions $\Phi^i$ for
\thetag{1.3} are determined implicitly through Lagrange function $L$:
$$
\sum^n_{s=1}\frac{\partial^2L}{\partial v^i\,\partial x^s}\,v^s+
\sum^n_{s=1}\frac{\partial^2L}{\partial v^i\,\partial v^s}\,\Phi^s=
\frac{\partial L}{\partial x^i}.
$$
However, not all Newtonian dynamical systems \thetag{1.1} can be
represented in Lagrangian form \thetag{1.3}. Therefore in \cite{26}
we considered more general fiber-preserving maps \thetag{1.4}.
We keep symbols $\lambda$ and $\lambda^{-1}$ for them. First is given
as follows:
$$
\hskip -2em
\cases p_1=L_1(x^1,\ldots,x^n,v^1,\ldots,v^n),\\
.\ .\ .\ .\ .\ .\ .\ .\ .\ .\ .\ .\ .\ .\ .\ .\
.\ .\ .\ .\ .\ .\ .\ \\
p_n=L_n(x^1,\ldots,x^n,v^1,\ldots,v^n).\\
\endcases
\tag1.6
$$
Fiber preserving map $\lambda\!:TM\to T^*\!M$ given by functions
\thetag{1.6} is called {\bf generalized Legendre transformation}.
Inverse map is given by functions
$$
\hskip -2em
\cases v^1=V^1(x^1,\ldots,x^n,p_1,\ldots,p_n),\\
.\ .\ .\ .\ .\ .\ .\ .\ .\ .\ .\ .\ .\ .\ .\ .\
.\ .\ .\ .\ .\ .\ .\ \\
v^n=V^n(x^1,\ldots,x^n,p_1,\ldots,p_n).\\
\endcases
\tag1.7
$$
For the sake of simplicity we shall assume generalized Legendre maps
\thetag{1.4} given by functions \thetag{1.6} and \thetag{1.7} in local
chart to be establishing diffeomorphism of $TM$ and $T^*\!M$ such that
zero maps to zero in each fiber of these two bundles.\par
    Now let's treat \thetag{1.7} as change of variables. Substituting
\thetag{1.7} into differential equations \thetag{1.1}, we can transform
them to the following ones:
$$
\xalignat 2
&\hskip -2em
\dot x^i=V^i
&&\dot p_i=\Theta_i.
\tag1.8
\endxalignat
$$
Here functions $V^i=V^i(x^1,\ldots,x^n,p_1,\ldots,p_n)$ are given by
\thetag{1.7}, while functions $\Theta^i=\Theta^i(x^1,\ldots,x^n,p_1,
\ldots,p_n)$ form another set of $n$ functions. They are determined
implicitly through function $\Phi^1,\,\ldots,\,\Phi^n$ in \thetag{1.1}
by the following equation:
$$
\hskip -2em
\Phi^i\compos\,\lambda^{-1}=\sum^n_{s=1}\frac{\partial V^i}{\partial x^s}
\,V^s+\sum^n_{s=1}\frac{\partial V^i}{\partial p_s}\,\Theta^s.
\tag1.9
$$
If the equations \thetag{1.9} are fulfilled, then both \thetag{1.1}
and \thetag{1.8} express the same dynamics, but in different
representations: \thetag{1.8} is called {\bf momentum representation}
or {\bf $\bold p$-representation} for Newtonian dynamical system
\thetag{1.1}, while \thetag{1.1} is called {\bf velocity representation}
or {\bf $\bold v$-representation} for Newtonian dynamical system
\thetag{1.8}. Concluding this section, note that \thetag{1.9} can be
written as
$$
\hskip -2em
\Theta_i\compos\lambda=\sum^n_{s=1}\frac{\partial L_i}{\partial x^s}
\,v^s+\sum^n_{s=1}\frac{\partial L_i}{\partial v^s}\,\Phi^s.
\tag1.10
$$
Formula \thetag{1.10} is more explicit with respect to functions
$\Theta_1,\,\ldots,\,\Theta_n$ in \thetag{1.8}.
\head
2. Extended tensor fields.
\endhead
    Extended tensor fields are closely related to generalized Legendre
maps \thetag{1.4}. Indeed, functions \thetag{1.6} are components
covector $\bold L\i T^*_p(M)$. However, unlike to components of
traditional covector field, they depend not only on coordinates
of point $p\in M$, but also on components of velocity vector
$\bold v\in T_p(M)$. Pair $q=(p,\bold v)$ is a point of $TM$.
Therefore $\bold L$ is extended covector field in
$\bold v$-representation (see definition below). In a similar
way, functions \thetag{1.7} are components of extended vector
field $\bold V$ in $\bold p$-representation.\par
    At each point $p$ of base manifold $M$ one can define tensor
space $T^r_s(p,M)$. This is the following tensor product of several
copies of $T_p(M)$ and $T^*_p(M)$:
$$
\hskip -2em
T^r_s(p,M)=\overbrace{T_p(M)\otimes\ldots
\otimes T_p(M)}^{\text{$r$ times}}\otimes
\underbrace{T^*_p(M)\otimes\ldots
\otimes T^*_p(M)}_{\text{$s$ times}}.
\tag2.1
$$
\definition{Definition 2.1} Extended tensor field $\bold X$
of type $(r,s)$ in $\bold v$-representation is a tensor-valued
function $\bold X=\bold X(q)$ with argument $q=(p,\bold v)$ in
tangent bundle $TM$ and with values in tensor space \thetag{2.1},
where $p=\pi(q)$.
\enddefinition
\definition{Definition 2.2} Extended tensor field $\bold X$
of type $(r,s)$ in $\bold p$-representation is a tensor-valued
function $\bold X=\bold X(q)$ with argument $q=(p,\bold p)$ in
cotangent bundle $T^*\!M$ and with values in tensor space
\thetag{2.1}, \pagebreak where $p=\pi(q)$.
\enddefinition
    Each extended tensor field can be transformed from $\bold v$
to $\bold p$-representation and back from $\bold p$ to
$\bold v$-representation by changing its argument. We use
generalized Legendre maps \thetag{1.4} for this purpose.
Extended covector field $\bold L$ and extended vector field
$\bold V$ defining these two maps are also examples of extended
tensor fields in $\bold v$ and $\bold p$-representations
respectively. Like traditional tensor fields, extended tensor
field $\bold X$ in local chart are represented by its components:
$$
\bold X=\sum^n_{i_1=1}\!...\!\sum^n_{i_r\,=1}\sum^n_{j_1=1}\!...\!
\sum^n_{j_s\,=1}X^{i_1\ldots\,i_r}_{j_1\ldots\,j_s}\,\frac{\partial}
{\partial x^{i_1}}\otimes\ldots\otimes\frac{\partial}{\partial x^{i_r}}
\otimes dx^{j_1}\otimes\ldots\otimes dx^{j_s}.
$$
\definition{Definition 2.3} Extended tensor field $\bold X$ is called
smooth if its components in any local chart of the manifold $M$ are
smooth functions of their arguments.
\enddefinition
Let's denote by $T^r_s(M)$ the set of smooth extended tensor fields of
type $(r,s)$ either in $\bold v$ or in $\bold p$-representation. They
form a module over the ring of smooth scalar functions in $TM$ or in
$T^*\!M$, we denote these rings by $\goth F(TM)$ and $\goth F(T^*\!M)$
respectively. Then we can define the following direct sum:
$$
\hskip -2em
\bold T(M)=\bigoplus^\infty_{r=0}\bigoplus^\infty_{s=0}T^r_s(M).
\tag2.2
$$
This direct sum \thetag{2.2} is a graded algebra over one of these two
rings. It is called {\bf extended algebra of tensor fields}. In $\bold
T(M)$ we have all standard tensorial operations like summation,
multiplication by scalars, tensor product, and contraction. Apart from
these algebraic operations, in $\bold T(M)$ we have {\bf canonical
vertical differentiation $\tilde\nabla$}. In $\bold v$-representation
it can be defined by the following formula:
$$
\hskip -2em
\tilde\nabla_{\!k}X^{i_1\ldots\,i_r}_{j_1\ldots\,j_s}=
\frac{\partial X^{i_1\ldots\,i_r}_{j_1\ldots\,j_s}}
{\partial v^k}.
\tag2.3
$$
For extended tensor field $\bold X$ in $\bold p$-representation we have
formula similar to \thetag{2.3}:
$$
\hskip -2em
\tilde\nabla^k\!X^{i_1\ldots\,i_r}_{j_1\ldots\,j_s}=
\frac{\partial X^{i_1\ldots\,i_r}_{j_1\ldots\,j_s}}
{\partial p_k}.
\tag2.4
$$
Canonical vertical differentiation does not commute with generalized
legendre transformation $\lambda$. Here we have the following
relationships:
$$
\align
&\hskip -2em
\tilde\nabla_{\!k}X^{i_1\ldots\,i_r}_{j_1\ldots\,j_s}
=\sum^n_{q=1}g_{qk}\cdot\tilde\nabla^q\!\!\left(X^{i_1\ldots\,i_r}_{j_1
\ldots\,j_s}\compos\lambda^{-1}\right)\!\compos\lambda,
\tag2.5
\\
&\hskip -2em
\tilde\nabla^k\!X^{i_1\ldots\,i_r}_{j_1\ldots\,j_s}
=\sum^n_{q=1}(g^{qk}\compos\lambda^{-1})\cdot\tilde\nabla_{\!q}
\kern-2.5pt\left(X^{i_1\ldots\,i_r}_{j_1\ldots\,j_s}\compos
\lambda\right)\!\compos\lambda^{-1}.
\tag2.6
\endalign
$$
By $g_{qk}$ and $g^{qk}$ in \thetag{2.5} and \thetag{2.6} we denote
Jacobi matrices for maps \thetag{1.4}:
$$
\pagebreak
\xalignat 2
&\hskip -2em
g_{qk}=\tilde\nabla_{\!k}L_q,
&&g^{qk}=\tilde\nabla^kV^q\compos\lambda.
\tag2.7
\endxalignat
$$
Matrices \thetag{2.7} are inverse to each other in the sense of the
following equalities:
$$
\xalignat 2
&\hskip -2em
\sum^n_{r=1}g^{sr}\,g_{rk}=\delta^s_k,
&&\sum^n_{r=1}g_{kr}\,g^{rs}=\delta^s_k.
\tag2.8
\endxalignat
$$
Extended tensor fields \thetag{2.7} are non-symmetric. But nevertheless,
looking to \thetag{2.5}, \thetag{2.6}, and \thetag{2.8}, we conclude
that they are analogs of metric tensor and dual metric tensor in
Riemannian geometry.
\head
3. Extended connections\\
and horizontal covariant derivatives.
\endhead
    Normality equations we are going to study below are written in
terms of covariant derivatives (see \cite{26}). Apart from \thetag{2.4},
there another covariant derivative is used. When applied to extended
tensor field $\bold X$ in $\bold p$-representation, it acts as follows:
$$
\hskip -2em
\aligned
&\nabla_{\!m}X^{i_1\ldots\,i_r}_{j_1\ldots\,j_s}=\frac{\partial
X^{i_1\ldots\,i_r}_{j_1\ldots\,j_s}}{\partial x^m}
+\sum^n_{a=1}\sum^n_{b=1}p_a\,\Gamma^a_{mb}\,\frac{\partial
X^{i_1\ldots\,i_r}_{j_1\ldots\,j_s}}{\partial p_b}\,+\\
&+\sum^r_{k=1}\sum^n_{a_k=1}\!\Gamma^{i_k}_{m\,a_k}\,X^{i_1\ldots\,
a_k\ldots\,i_r}_{j_1\ldots\,\ldots\,\ldots\,j_s}
-\sum^s_{k=1}\sum^n_{b_k=1}\!\Gamma^{b_k}_{m\,j_k}
X^{i_1\ldots\,\ldots\,\ldots\,i_r}_{j_1\ldots\,b_k\ldots\,j_s}.
\endaligned
\tag3.1
$$
Here $\Gamma^k_{ij}$ are components of some extended affine
connection\footnotemark\ in $\bold p$-representation:
$$
\hskip -2em
\Gamma^k_{ij}=\Gamma^k_{ij}(x^1,\ldots,x^n,p_1,\ldots,p_n).
\tag3.2
$$
\footnotetext{Below we consider only symmetric connections
$\Gamma^k_{ij}=\Gamma^k_{ji}$. This is sufficient for our purposes.}
\adjustfootnotemark{-1}
Covariant differentiation $\nabla$ in $\bold T(M)$ determined by
extended affine connection \thetag{3.2} in formula \thetag{3.1}
is called {\bf horizontal covariant differentiation}. It is not
canonical since it depends on the choice of $\Gamma$. Now we shall
not discuss the concept of extended affine connection (see \cite{26}
and Chapter~\uppercase\expandafter{\romannumeral 3} of thesis
\cite{15}). Note only that by means of generalized Legendre map
$\lambda$ we can transform its components \thetag{3.2} to
$\bold v$-representation: $\Gamma^k_{ij}\to\Gamma^k_{ij}\compos
\lambda$. Then we have
$$
\hskip -2em
\Gamma^k_{ij}=\Gamma^k_{ij}(x^1,\ldots,x^n,v^1,\ldots,v^n).
\tag3.3
$$
Using \thetag{3.3}, one can define horizontal covariant differentiation
in $\bold v$-representation:
$$
\hskip -2em
\aligned
&\nabla_{\!m}X^{i_1\ldots\,i_r}_{j_1\ldots\,j_s}=\frac{\partial
X^{i_1\ldots\,i_r}_{j_1\ldots\,j_s}}{\partial x^m}
-\sum^n_{a=1}\sum^n_{b=1}v^a\,\Gamma^b_{am}\,\frac{\partial
X^{i_1\ldots\,i_r}_{j_1\ldots\,j_s}}{\partial v^b}\,+\\
&+\sum^r_{k=1}\sum^n_{a_k=1}\!\Gamma^{i_k}_{m\,a_k}\,X^{i_1\ldots\,
a_k\ldots\,i_r}_{j_1\ldots\,\ldots\,\ldots\,j_s}
-\sum^s_{k=1}\sum^n_{b_k=1}\!\Gamma^{b_k}_{m\,j_k}
X^{i_1\ldots\,\ldots\,\ldots\,i_r}_{j_1\ldots\,b_k\ldots\,j_s}.
\endaligned
\tag3.4
$$
Horizontal covariant differentiation $\nabla$ also does not commute
\pagebreak with generalized Legendre map $\lambda$. Here for
differentiation $\nabla$ we have the following relationships:
$$
\align
&\hskip -2em
\aligned
\nabla_{\!m}X^{i_1\ldots\,i_r}_{j_1\ldots\,j_s}&=
\nabla_{\!m}\kern-2.5pt\left(X^{i_1\ldots\,i_r}_{j_1\ldots\,j_s}
\compos\lambda\right)\!\compos\lambda^{-1}+\\
&+\sum^n_{q=1}\nabla_{\!m}V^q\cdot\tilde\nabla_{\!q}\kern-2.5pt
\left(X^{i_1\ldots\,i_r}_{j_1\ldots\,j_s}\compos\lambda\right)\!\compos
\lambda^{-1},
\endaligned
\tag3.5
\\
&\hskip -2em
\aligned
\nabla_{\!m}X^{i_1\ldots\,i_r}_{j_1\ldots\,j_s}&=
\nabla_{\!m}\kern-2.5pt\left(X^{i_1\ldots\,i_r}_{j_1\ldots\,j_s}
\compos\lambda^{-1}\right)\!\compos\lambda+\\
&+\sum^n_{q=1}\nabla_{\!m}L_q\cdot\tilde\nabla^q\kern-2.5pt
\left(X^{i_1\ldots\,i_r}_{j_1\ldots\,j_s}\compos\lambda^{-1}\right)
\!\compos\lambda.
\endaligned
\tag3.6
\endalign
$$
These relationships \thetag{3.5} and \thetag{3.6} for $\nabla$ are
analogs of relationships \thetag{2.5} and \thetag{2.6}.
\head
4. Force vector and force covector
of Newtonian dynamical system.
\endhead
    Note that functions $\Phi^i(x^1,\ldots,x^n,v^1,\ldots,v^n)$ in
\thetag{1.1} are not components of vector field even if we understand
it in the sense of definition~2.1. Similarly, functions $\Theta_i(x^1,
\ldots,x^n,p_1,\ldots,p_n)$ in \thetag{1.8} are not components of
covector field in the sense of definition~2.2. We can change this
situation if we choose some symmetric extended affine connection
$\Gamma$. Using components of $\Gamma$ in $\bold v$-representation
\thetag{3.3}, we replace time derivative $\dot v^i$ in \thetag{1.1}
by covariant time derivative
$$
\hskip -2em
\nabla_tv^i=\dot v^i+\sum^n_{j=1}\sum^n_{k=1}\Gamma^i_{jk}\,v^j\,v^k.
\tag4.1
$$
In a similar way we replace time derivative $\dot p_i$ in \thetag{1.8}
by covariant time derivative
$$
\hskip -2em
\nabla_tp_i=\dot p_i-\sum^n_{j=1}\sum^n_{k=1}\Gamma^k_{ij}\,V^j\,p_k.
\tag4.2
$$
For this purpose we need components of $\Gamma$ in
$\bold p$-representation \thetag{3.2}. Relying upon \thetag{4.1} and
\thetag{4.2} we replace functions $\Phi^i$ and $\Theta_i$ by functions
$F^i$ and $Q_i$ as follows:
$$
\align
&\hskip -2em
F^i=\Phi^i+\sum^n_{j=1}\sum^n_{k=1}\Gamma^i_{jk}\,v^j\,v^k,
\tag4.3\\
&\hskip -2em
Q_i=\Theta_i-\sum^n_{j=1}\sum^n_{k=1}\Gamma^k_{ij}\,V^j\,p_k.
\tag4.4
\endalign
$$
In terms or these newly introduced functions \thetag{4.3}, \thetag{4.4}
and covariant derivatives \thetag{4.1}, \thetag{4.2} differential
equations \thetag{1.1} and \thetag{1.8} are rewritten as
$$
\xalignat 2
&\hskip -2em
\dot x^i=v^i,
&&\nabla_tv^i=F^i,
\tag4.5\\
&\hskip -2em
\dot x^i=V^i,
&&\nabla_tp_i=Q_i.
\tag4.6
\endxalignat
$$
We say that \thetag{4.5} and \thetag{4.6} \pagebreak are tensorial form of
differential equations of Newtonian dynamics in $\bold v$-representation
and $\bold p$-representation respectively. Indeed, functions $F^i=F^i(x^1,
\ldots,x^n,v^1,\ldots,v^n)$ in \thetag{4.5} are components of extended
vector field $\bold F$. It is called {\bf force vector} of Newtonian
dynamical system \thetag{4.5}. Functions $Q_i=Q_i(x^1,\ldots,x^n,v^1,
\ldots,v^n)$ represent extended covector field $\bold Q$. It is called
{\bf force covector} of Newtonian dynamical system \thetag{4.6}.
\head
5. Normal shift and normality equations.
\endhead
    In Euclidean and in Riemannian geometry normal shift is a continuous
deformation of hypersurface when each its point moves along some
trajectory and moving hypersurface in whole keeps orthogonality to
trajectories of all its points. This means that normal vector of moving
hypersurface is always collinear to tangent vector of shift trajectory.
The same is true in Finslerian geometry (see
Chapter~\uppercase\expandafter{\romannumeral 8} of thesis \cite{15}).
In non-metric geometries one cannot define normal vector of hypersurface.
But here one can define {\bf normal covector}. This is key idea of
papers \cite{24}, \cite{25}, and \cite{26}. Let $\sigma$ be some
hypersurface in $M$ and let $p\in\sigma$ be a point of $\sigma$. Then
$T_p(\sigma)$ is a hyperplane in tangent space $T_p(M)$. Denote by
$\bold n$ any nonzero covector in $T^*_p(M)$ whose null-space coincides
with $T_p(\sigma)$:
$$
\hskip -2em
T_p(\sigma)=\{\bold X\in T_p(P): \left<\bold n\,|\,\bold X\right>=0\}.
\tag5.1
$$
By angular brackets in \thetag{5.1} we denote scalar product of vector
$\bold X$ and covector $\bold n$:
$$
\hskip -2em
\left<\bold n\,|\,\bold X\right>=\sum^n_{i=1}n_i\,X^i.
\tag5.2
$$
Covector $\bold n=\bold n(p)$ is determined by the condition \thetag{5.1}
uniquely up to some nonzero scalar factor: $\bold n\to\nu\cdot\bold n$.
It is called {\bf normal covector} of hypersurface $\sigma$ at the point
$p$. One can choose $\bold n=\bold n(p)$ to be smooth covector function
on $\sigma$ (at least locally in some neighborhood of any point).
Now we can associate with $\sigma$ the following initial data for
Newtonian dynamical system \thetag{1.8} in $\bold p$-representation:
$$
\xalignat 2
&\hskip -2em
x^i\,\hbox{\vrule height 8pt depth 8pt width 0.5pt}_{\,t=0}
=x^i(p),
&&p_i\,\hbox{\vrule height 8pt depth 8pt width 0.5pt}_{\,t=0}
=\nu(p)\cdot n_i(p).
\tag5.3
\endxalignat
$$
Let's fix some point $p_0\in\sigma$. Smooth normal covector $\bold n
=\bold n(p)$ of hypersurface $\sigma$ is defined in some neighborhood
of the point $p_0$. For scalar factor $\nu=\nu(p)$ in \thetag{5.3}
we shall assume it to be smooth nonzero function in some possibly
smaller neighborhood of the point $p_0$ normalized by the condition
$$
\hskip -2em
\nu(p_0)=\nu_0,
\tag5.4
$$
where $\nu_0\neq 0$ is some arbitrary predefined constant. We can apply
initial data \thetag{5.3} to Newtonian dynamical system given either
by differential equations \thetag{1.8} or \thetag{4.6}. This initiates
motion of hypersurface $\sigma$ (or its smaller part including our fixed
point $p_0$) when all points are moving along trajectories of this
dynamical system.
\definition{Definition 5.1} Shift of hypersurface $\sigma\subset M$
along trajectories of a dynamical system in cotangent bundle $T^*\!M$
determined by initial data \thetag{5.3} is called {\bf normal shift}
if normal covectors of all shifted hypersurfaces $\sigma_t$ are
collinear to momentum covector $\bold p$ \pagebreak on shift
trajectories.
\enddefinition
\definition{Definition 5.2} Dynamical system in cotangent bundle
$T^*\!M$ satisfies {\bf strong normality condition} if for any
hypersurface $\sigma\subset M$, for any fixed point $p_0\in\sigma$,
and for any constant $\nu_0\neq 0$ there is some open neighborhood
of $p_0$ on $\sigma$ and some smooth function $\nu(p)$ normalized
by the condition \thetag{5.4} in this neighborhood such that it
determines normal shift of $\sigma$ in the sense of definition~5.1.
\enddefinition
   We say that dynamical system {admits normal shift of hypersurfaces}
in $M$ if it satisfies strong normality condition stated in
definition~5.2. It is clear that when applied to Newtonian dynamical
system \thetag{4.6}, strong normality condition yields some restrictions
for the functions $V^1,\,\ldots,\,V^n$ and $Q_1,\,\ldots,\,Q_n$ in
\thetag{4.6}. In \cite{26} it was shown that these restrictions are
written in form of so called {\bf normality equations} for these
functions. In order to write them we need to introduce several auxiliary
extended tensor fields. First is a vector field $\bold W$ with
components
$$
\hskip -2em
W^i=\sum^n_{s=1}p_s\,\tilde\nabla^iV^s.
\tag5.5
$$
This vector field is used to define projector-valued operator field
$\bold P$ with components
$$
\hskip -2em
P^i_j=\delta^i_j-\frac{W^i\,p_j}{\Omega},
\tag5.6
$$
where $\Omega$ is a scalar field defined as scalar product of covector
$\bold p$ and vector $\bold W$:
$$
\hskip -2em
\Omega=\left<\bold p\,|\,\bold W\right>=\sum^n_{s=1}p_s\,W^s
\tag5.7
$$
(compare with \thetag{5.2}). Next is covector field $\bold U$.
Its components are given by formula
$$
\hskip -2em
U_i=Q_i+\sum^n_{s=1}\nabla_{\!i}V^s\,p_s.
\tag5.8
$$
Formulas for vector field $\boldsymbol\alpha$ and covector field
$\boldsymbol\beta$ are more complicated:
$$
\gather
\hskip -2em
\gathered
\alpha^k=\sum^n_{r=1}\tilde\nabla^kV^r\,U_r+\sum^n_{r=1}\nabla_{\!r}
W^k\,V^r+\sum^n_{r=1}\tilde\nabla^rW^k\,Q_r\,+\\
+\,\sum^n_{r=1}W^r\,\tilde\nabla^kQ_r-\sum^n_{r=1}\sum^n_{s=1}
\sum^n_{q=1}p_s\,D^{sk}_{rq}\,W^r\,V^q.
\endgathered
\tag5.9
\\
\hskip -2em
\gathered
\beta_k=\sum^n_{r=1}\nabla_{\!r}U_k\,V^r+\sum^n_{r=1}
\tilde\nabla^rU_k\,Q_r+\sum^n_{r=1}\nabla_{\!k}V^r\,U_r\,+\\
+\sum^n_{r=1}\nabla_{\!k}Q_r\,W^r-\sum^n_{r=1}\sum^n_{s=1}
\sum^n_{m=1}\left(R^s_{rmk}\,V^m-D^{sm}_{rk}\,Q_m\right)W^r\,p_s.
\endgathered
\tag5.10
\endgather
$$
Here in \thetag{5.9} and \thetag{5.10} we see components of two
curvature tensors $\bold D$ and $\bold R$:
$$
\hskip -2em
D^{kr}_{ij}=-\frac{\partial\Gamma^k_{ij}}{\partial p_r}.
\tag5.11
$$
Tensor $\bold D$ is dynamic curvature tensor. It is nonzero only for
extended connections. As for tensor $\bold R$, it is analog of standard
curvature tensor:
$$
\gathered
R^k_{rij}=\frac{\partial\Gamma^k_{jr}}{\partial x^i}-
\frac{\partial\Gamma^k_{ir}}{\partial x^j}+\sum^n_{m=1}
\Gamma^k_{im}\,\Gamma^m_{jr}-\sum^n_{m=1}\Gamma^k_{jm}\,
\Gamma^m_{ir}\,+\\
\vspace{1ex}
+\sum^n_{m=1}\sum^n_{s=1}p_s\,\Gamma^s_{mi}\,
\frac{\partial\Gamma^k_{jr}}{\partial p^m}-\sum^n_{m=1}
\sum^n_{s=1}p_s\,\Gamma^s_{mj}\,\frac{\partial
\Gamma^k_{ir}}{\partial p^m}.
\endgathered
\tag5.12
$$
Covector field $\boldsymbol\eta$ is defined by the following formula
for its components:
$$
\hskip -2em
\eta_k=\beta_k-\sum^n_{s=1}\frac{U_k\,\alpha^s\,p_s}{\Omega}.
\tag5.13
$$\par
As it was shown in \cite{26} in multidimensional case $n=\dim M
\geqslant 3$ complete system of normality equations naturally
subdivides into two parts. First part is formed by so called
{\bf weak normality equations}. Relying upon all the above
notations \thetag{5.5}, \thetag{5.6}, \thetag{5.7}, \thetag{5.8},
\thetag{5.9}, \thetag{5.10}, \thetag{5.11}, \thetag{5.12}, and
\thetag{5.13}, now we are able to write weak normality equations
in a very concise form:
$$
\xalignat 2
\hskip -2em
&\sum^n_{r=1}\alpha^r\,P^k_r=0,
&&\sum^n_{r=1}\eta_r\,P^r_k=0.
\tag5.14
\endxalignat
$$
Inspecting \thetag{5.14} with all the above notations in mind, one
can find that weak normality equations form a system of $2n$ partial
differential equations for components of two fields: extended vector
field $\bold V$ and extended covector field $\bold Q$. But rank of
projection operator $\bold P$ is $n-1$. Therefore actual number of
independent equations in \thetag{5.14} is equal to $2n-2$.\par
    Second part of normality equations in multidimensional case
$n=\dim M\geqslant 3$ is formed by so called {\bf additional normality
equations}. In order to write them we need to use some more auxiliary
notations:
$$
\align
&\hskip -2em
A^{rs}=\tilde\nabla^rW^s,
\tag5.15\\
&\hskip -2em
\aligned
B^r_s&=\tilde\nabla^rU_s+\sum^n_{m=1}\sum^n_{k=1}W^k\,p_m\,D^{mr}_{ks}\,-\\
&-\,\nabla_{\!s}W^r+\sum^n_{m=1}\frac{\tilde\nabla^mW^r-\tilde\nabla^rW^m}
{\Omega}\,U_s\,p_m,
\endaligned
\tag5.16\\
&\hskip -2em
\aligned
C_{rs}&=\nabla_{\!r}U_s
-\sum^n_{m=1}\frac{U_r\,\tilde\nabla^mU_s+U_s\,\nabla_{\!r}W^m}
{\Omega}\,p_m\,-
\\
&-\sum^n_{k=1}\sum^n_{q=1}\left(\,\shave{\sum^n_{m=1}}\frac{D^{mq}_{ks}
\,U_r}{\Omega}\,p_m+\frac{R^q_{krs}}{2}\right)W^k\,p_q.
\endaligned
\tag5.17
\endalign
$$
Here $A^{rs}$, $B^r_s$, and $C_{rs}$ are components of three
extended tensor fields $\bold A$, $\bold B$, and $\bold C$
respectively. Now, using \thetag{5.15}, \thetag{5.16}, \thetag{5.17}
and keeping in mind all previous notations, we can write additional
normality equations as well. They look like
$$
\align
&\hskip -2em
\sum^n_{r=1}\sum^n_{s=1}(A^{rs}-A^{sr})\,P^i_r\,P^j_s=0,
\tag5.18\\
&\hskip -2em
\sum^n_{r=1}\sum^n_{s=1}P^i_r\,B^r_s\,P^s_j=
\lambda\,P^i_j,
\tag5.19\\
&\hskip -2em
\sum^n_{r=1}\sum^n_{s=1}(C_{rs}-C_{sr})\,P^r_i\,P^s_j=0.
\tag5.20
\endalign
$$
Here $\lambda$ in \thetag{5.19} is a scalar factor. However, it
is not a parameter with undetermined or deliberate value. Its
value is determined by the equation \thetag{5.19} itself:
$$
\lambda=\frac{\tr(\lambda\,\bold P)}{\tr\bold P}=\sum^n_{r=1}
\sum^n_{s=1}\frac{B^r_s\,P^s_r}{n-1}.
\tag5.21
$$
Due to formula \thetag{5.21} the equation \thetag{5.19} can be
written as follows:
$$
\hskip -2em
\sum^n_{r=1}\sum^n_{s=1}P^i_r\,B^r_s\,P^s_j=\sum^n_{r=1}\sum^n_{s=1}
\frac{B^r_s\,P^s_r}{n-1}\,P^i_j.
\tag5.22
$$
Under some additional restriction for generalized Legendre map
(see definition~6.1 in paper \cite{26}) the role of weak and
additional normality equations is described by the following
result (see theorems~8.1, 11.1, and 12.1 in paper \cite{26}).
\proclaim{Theorem 5.1} Strong normality condition for Newtonian
dynamical system \thetag{4.6} in multidimensional case $n\geqslant
3$ is equivalent to complete system of normality equations
\thetag{5.14}, \thetag{5.18}, \thetag{5.19}, \thetag{5.20}
that should be fulfilled at all points $q=(p,\bold p)$ of
cotangent bundle $T^*\!M$, where $\bold p\neq 0$.
\endproclaim
Two-dimensional case is an exception. Here we have the following
result.
\proclaim{Theorem 5.2} Strong normality condition for Newtonian
dynamical system \thetag{4.6} in the dimension $n=2$ is equivalent
to weak normality equations \thetag{5.14} that should be fulfilled
at all points $q=(p,\bold p)$ of cotangent bundle $T^*\!M$, where
$\bold p\neq 0$.
\endproclaim
    Weak normality equations \thetag{5.14} and additional normality
equations \thetag{5.18}, \thetag{5.19}, \thetag{5.20}, as well as
the equations \thetag{5.22}, are written in terms of extended tensor
fields and covariant derivatives in $\bold p$-representation. Our
goal below is to transform all of them to $\bold v$-representation.
\head
6. Transformation of projection operator.
\endhead
    Projection operator $\bold P$ with components \thetag{5.6}
is present in each normality equation. Let's transform it to
$\bold v$-representation. For $W^i$ in \thetag{5.6} we have
$$
\hskip -2em
W^i\compos\lambda=\sum^n_{s=1}(\,p_s\compos\lambda)\,
(\tilde\nabla^iV^s\compos\lambda)=\sum^n_{s=1}L_s\,g^{si}=L^i.
\tag6.1
$$
Here we used formula \thetag{5.5} for $W^i$ and notations \thetag{2.7}.
Formula $p_s\compos\lambda=L_s$ is quite obvious (see formulas \thetag{1.6},
where generalized Legendre map $\lambda$ is defined by functions $L_1,\,
\ldots,\,L_n$ in local chart). By $L^i$ in \thetag{6.1} we denote components
of extended vector field $\bold L$ dual to extended covector field $\bold L$
with respect to non-symmetric extended metric tensor \thetag{2.7}.\par
    Now let's transform extended scalar field $\Omega$ that forms
denominator of fraction in formula \thetag{5.6}. Applying $\lambda$ to
formula \thetag{5.7}, we derive
$$
\hskip -2em
\Omega\,\compos\lambda=\sum^n_{s=1}(\,p_s\compos\lambda)\,(W^s\compos\lambda)
=\sum^n_{s=1}L_s\,L^s=|\bold L|^2.
\tag6.2
$$
For the last sum in \thetag{6.2} by similarity we used the same notation
$|\bold L|^2$ as in Riemannian geometry, though the length of vector $\bold
L$ here is understood in the sense of extended metric \thetag{2.7}.
Combining \thetag{6.1} and \thetag{6.2}, for components of $\bold P$ we get
$$
\hskip -2em
P^i_j\compos\lambda=\delta^i_j-\frac{L^i\,L_j}{|\bold L|^2}.
\tag6.3
$$
Note that formula \thetag{6.3} for $P^i_j$ looks like formula for
components of orthogonal projector in Riemannian geometry.
Writing $\Omega\,\compos\lambda$ and $P^i_j\compos\lambda$ in \thetag{6.2}
and \thetag{6.3} means that we take $\bold v$-representation of extended
tensor fields $\Omega$ and $\bold P$, which are initially in
$\bold p$-representation. But having transformed them to
$\bold v$-representation, we can then omit $\lambda$ symbol, i\.\,e\. we
can write
$$
\xalignat 2
&\hskip -2em
\Omega=\sum^n_{s=1}L_s\,L^s=|\bold L|^2,
&&P^i_j=\delta^i_j-\frac{L^i\,L_j}{|\bold L|^2},
\tag6.4
\endxalignat
$$
assuming both sides of the equalities \thetag{6.4} to be in
$\bold v$-representation. Below we shall use the same convention
for other extended tensor fields $\bold U$, $\boldsymbol\alpha$,
$\boldsymbol\beta$, $\boldsymbol\eta$, $\bold A$, $\bold B$, and
$\bold C$.
\head
7. Weak normality equations in $\bold v$-representation.
\endhead
    Now let's continue transforming tensor fields used in writing
weak normality equations \thetag{5.14}. For $\nabla_{\!i}V^s$ in
\thetag{5.8}, applying formula \thetag{3.6}, we derive
$$
\hskip -2em
\nabla_{\!i}v^s=\nabla_{\!i}V^s\compos\lambda+\sum^n_{q=1}
\nabla_{\!i}L_q\,(\tilde\nabla^qV^s\compos\lambda)
=\nabla_{\!i}V^s\compos\lambda+\sum^n_{q=1}\nabla_{\!i}L_q\ g^{sq}.
\tag7.1
$$
Left hand side of \thetag{7.1} is equal to zero. Therefore we have
the equality
$$
\nabla_{\!i}V^s\compos\lambda=-\sum^n_{q=1}\nabla_{\!i}L_q\ g^{sq}.
\tag7.2
$$
Now we can substitute \thetag{7.2} into \thetag{5.8}. Taking into
account \thetag{6.1}, then we get
$$
\hskip -2em
U_i\compos\lambda=Q_i\compos\lambda-\sum^n_{q=1}L^q\,\nabla_{\!i}L_q.
\tag7.3
$$
Components of covector field $\bold Q$ in $\bold p$-representation
define Newtonian dynamical system \thetag{4.6}. However, its components
in $\bold v$-representation have no direct relation to \thetag{4.5}.
This means that $Q_i\compos\lambda$ in \thetag{7.3} should be expressed
through $F^i$. Combining \thetag{1.10} with formulas \thetag{4.3} and
\thetag{4.4}, we get the following expression for $Q_i\compos\lambda$:
$$
\hskip -2em
Q_i\compos\lambda=\sum^n_{q=1}v^q\,\nabla_{\!q}L_i+
\sum^n_{q=1}F^q\,\tilde\nabla_{\!q}L_i=\sum^n_{q=1}
v^q\,\nabla_{\!q}L_i+\sum^n_{q=1}F^q\,g_{iq}.
\tag7.4
$$
By introducing extended covector field $\bold F$ dual to extended
vector field $\bold F$ with respect to metric \thetag{2.7} we can
simplify \thetag{7.4} a little bit more
$$
\hskip -2em
Q_i\compos\lambda=\sum^n_{q=1}v^q\,\nabla_{\!q}L_i+F_i.
\tag7.5
$$
For $L_i$, $L^q$, $F_i$, and $F^q$ we have the following relationships:
$$
\xalignat 2
&\hskip -2em
L^i=\sum^n_{q=1}L_q\,g^{qi},
&&L_q=\sum^n_{i=1}L^i\,g_{iq},
\tag7.6\\
&\hskip -2em
F^i=\sum^n_{q=1}g^{iq}\,F_q,
&&F_q=\sum^n_{i=1}g_{qi}\,F^i.
\tag7.7
\endxalignat
$$
Comparing \thetag{7.6} and \thetag{7.7}, we see that due to asymmetry
of metric \thetag{2.7} they are slightly different. For the sake of
certainty we say that $L^i$ and $L_q$ are right-dual to each other,
while $F^i$ and $F_q$ are left-dual. Substituting \thetag{7.5} into
\thetag{7.3}, we get
$$
\hskip -2em
U_i=\sum^n_{q=1}v^q\,\nabla_{\!q}L_i
-\sum^n_{q=1}L^q\,\nabla_{\!i}L_q+F_i.
\tag7.8
$$
According to our convention, we omitted symbol $\lambda$ in left hand
side of \thetag{7.8} since this is ultimate expression for components
of covector field $\bold U$ in $\bold v$-representation.\par
    Next step is to transform extended vector field $\boldsymbol\alpha$
with components \thetag{5.9}. By means of formula \thetag{3.5} for
covariant derivative $\nabla_{\!r}W^k$ we derive:
$$
\hskip -2em
\nabla_{\!r}W^k\compos\lambda=\nabla_{\!r}L^k+\sum^n_{q=1}
(\nabla_{\!r}V^q\compos\lambda)\,\tilde\nabla_{\!q}L^k.
\tag7.9
$$
Now, if we apply formula \thetag{7.2} to $\nabla_{\!r}V^q\compos\lambda$
in formula \thetag{7.9}, we obtain
$$
\hskip -2em
\nabla_{\!r}W^k\compos\lambda=\nabla_{\!r}L^k-\sum^n_{q=1}
\sum^n_{s=1}\nabla_{\!r}L_s\,g^{qs}\,\tilde\nabla_{\!q}L^k.
\tag7.10
$$
Using formula \thetag{7.10}, we can transform second term in right hand
side of \thetag{5.9}:
$$
\sum^n_{r=1}(\nabla_{\!r}W^k\,V^r)\compos\lambda=\sum^n_{r=1}v^r\,
\nabla_{\!r}L^k-\sum^n_{q=1}\sum^n_{r=1}\sum^n_{s=1}v^r\,\nabla_{\!r}
L_s\,g^{qs}\,\tilde\nabla_{\!q}L^k.\quad
\tag7.11
$$
Using formulas \thetag{2.6} and \thetag{7.5}, for third term in right
hand side of \thetag{5.9} we derive
$$
\hskip -2em
\gathered
\sum^n_{r=1}(\tilde\nabla^rW^k\,Q_r)\compos\lambda=\sum^n_{q=1}
\sum^n_{r=1}g^{qr}\,\tilde\nabla_qL^k\,(Q_r\compos\lambda)=\\
=\sum^n_{q=1}\sum^n_{r=1}\sum^n_{s=1}g^{qr}\,\tilde\nabla_qL^k
\,v^s\,\nabla_{\!s}L_r+\sum^n_{q=1}F^q\,\tilde\nabla_qL^k.
\endgathered
\tag7.12
$$
Transforming first term in right hand side of \thetag{5.9} needs no
special efforts. Indeed, for $\tilde\nabla^kV^r$ we use formula
\thetag{2.7}. As a result we get
$$
\hskip -2em
\sum^n_{r=1}(\tilde\nabla^kV^r\,U_r)\compos\lambda=
\sum^n_{r=1}g^{rk}\,U_r.
\tag7.13
$$
Components of covector field $\bold U$ in right hand side of
\thetag{7.13} are already transformed to $\bold v$-representation,
see formula \thetag{7.8}.\par
    In fourth term in right hand side of \thetag{5.9} we see $\tilde
\nabla^kQ_r$. Applying formulas \thetag{2.6} and \thetag{7.5}, for
this derivative we obtain the following expression:
$$
\hskip -2em
\gathered
\tilde\nabla^kQ_r\compos\lambda=\sum^n_{q=1}g^{qk}\,\tilde\nabla_{\!q}
(Q_r\compos\lambda)=\sum^n_{q=1}g^{qk}\,\nabla_{\!q}L_r+\\
+\sum^n_{q=1}\sum^n_{s=1}g^{qk}\,v^s\,\tilde\nabla_{\!q}\nabla_{\!s}L_r
+\sum^n_{q=1}\sum^n_{s=1}g^{qk}\,\tilde\nabla_{\!q}F_r.
\endgathered
\tag7.14
$$
Using formula \thetag{7.14}, for the fourth term in right hand side
of \thetag{5.9} we derive
$$
\hskip -2em
\gathered
\sum^n_{r=1}(W^r\,\tilde\nabla^kQ_r)\compos\lambda=\sum^n_{r=1}
\sum^n_{q=1}L^r\,g^{qk}\,\nabla_{\!q}L_r+\\
+\sum^n_{q=1}\sum^n_{r=1}\sum^n_{s=1}L^r\,g^{qk}\,v^s\,
\tilde\nabla_{\!q}\nabla_{\!s}L_r
+\sum^n_{q=1}\sum^n_{r=1}\sum^n_{s=1}L^r\,g^{qk}\,\tilde
\nabla_{\!q}F_r.
\endgathered
\tag7.15
$$
Last term in right hand side of \thetag{5.9} contains components of
dynamic curvature tensor $\bold D$. In $\bold p$-representation they
are given by formula \thetag{5.11}. In $\bold v$-representation we
have analogous formula for components of dynamic curvature tensor:
$$
\hskip -2em
D^k_{rij}=-\frac{\partial\Gamma^k_{ir}}{\partial v^j}.
\tag7.16
$$
Note that \thetag{5.11} and \thetag{7.16} define two different
extended tensor fields (not two representations of the same field).
But there is very simple relationship binding components of these
two fields in $\bold p$ and $\bold v$-representations:
$$
\hskip -2em
D^{kr}_{ij}\compos\lambda=\sum^n_{s=1}g^{sr}\,D^k_{ijs}.
\tag7.17
$$
On the base of \thetag{7.17} for last term in right hand side of
\thetag{5.9} we derive:
$$
\sum^n_{r=1}\sum^n_{s=1}\sum^n_{q=1}(p_s\,D^{sk}_{rq}\,W^r\,V^q)
\compos\lambda=\sum^n_{q=1}\sum^n_{m=1}\sum^n_{r=1}\sum^n_{s=1}
g^{mk}\,L_s\,D^s_{rqm}\,L^r\,v^q.\quad
\tag7.18
$$
Now we can summarize above calculations \thetag{7.11}, \thetag{7.12},
\thetag{7.13}, \thetag{7.15}, \thetag{7.18} and write formula for
components of extended covector field $\boldsymbol\alpha$ in
$\bold v$-representation: 
$$
\hskip -2em
\gathered
\alpha^k=\sum^n_{r=1}g^{rk}\,U_r+\sum^n_{r=1}v^r\,\nabla_{\!r}L^k
+\sum^n_{q=1}F^q\,\tilde\nabla_qL^k\,+\\
+\sum^n_{q=1}\sum^n_{r=1}\sum^n_{s=1}L^r\,g^{qk}\,v^s\,
\tilde\nabla_{\!q}\nabla_{\!s}L_r
+\sum^n_{q=1}\sum^n_{r=1}L^r\,g^{qk}\,\tilde
\nabla_{\!q}F_r\,+\\
+\sum^n_{r=1}\sum^n_{q=1}L^r\,g^{qk}\,\nabla_{\!q}L_r
-\sum^n_{q=1}\sum^n_{m=1}\sum^n_{r=1}\sum^n_{s=1}
g^{mk}\,L_s\,D^s_{rqm}\,L^r\,v^q.
\endgathered
\tag7.19
$$\par
    Second normality equation in \thetag{5.14} is written in terms of
components of covector field $\boldsymbol\eta$. In order to use formula
\thetag{5.13} for $\eta_k$ we should first transform covector field
$\boldsymbol\beta$ with components \thetag{5.10} to
$\bold v$-representation. Using \thetag{3.5} and \thetag{7.2}, we get
$$
\hskip -2em
\nabla_{\!r}U_k\compos\lambda=\nabla_{\!r}U_k-\sum^n_{q=1}\sum^n_{s=1}
\nabla_{\!r}L_q\ g^{sq}\,\tilde\nabla_{\!s}U_k.
\tag7.20
$$
For $U_k$ in right hand side of \thetag{7.20} we should use
\thetag{7.8} since covector field $\boldsymbol U$ is already
transformed to $\bold v$-representation. Substituting \thetag{7.20}
into \thetag{5.10}, for the first term in right hand side of
\thetag{5.10} we derive
$$
\hskip -2em
\sum^n_{r=1}(\nabla_{\!r}U_k\,V^r)\compos\lambda=\sum^n_{r=1}
\nabla_{\!r}U_k\,v^r-\sum^n_{q=1}\sum^n_{r=1}\sum^n_{s=1}
\nabla_{\!r}L_q\ g^{sq}\,\tilde\nabla_{\!s}U_k\,v^r.
\tag7.21
$$
In order to transform $\tilde\nabla^rU_k$ to $\bold v$-representation
we apply formula \thetag{2.6}. This yields
$$
\hskip -2em
\tilde\nabla^rU_k\compos\lambda=\sum^n_{q=1}g^{qr}\,
\tilde\nabla_{\!q}U_k.
\tag7.22
$$
Applying \thetag{7.22} and \thetag{7.5}, for the second term in
right hand side of \thetag{5.10} we get
$$
\pagebreak
\hskip -2em
\gathered
\sum^n_{r=1}(\tilde\nabla^rU_k\,Q_r)\compos\lambda=\sum^n_{q=1}
\sum^n_{r=1}g^{qr}\,\tilde\nabla_{\!q}U_k\,F_r\,+\\
+\sum^n_{q=1}\sum^n_{r=1}\sum^n_{s=1}g^{qr}\,\tilde\nabla_{\!q}U_k
\,v^s\,\nabla_{\!s}L_r.
\endgathered
\tag7.23
$$
In order to transform third term in right hand side of \thetag{5.10}
we use \thetag{7.2} and \thetag{7.8}:
$$
\hskip -2em
\gathered
\sum^n_{r=1}\nabla_{\!k}V^r\,U_r=-\sum^n_{q=1}\sum^n_{r=1}\sum^n_{s=1}
\nabla_{\!k}L_q\,g^{rq}\,v^s\,\nabla_{\!s}L_r\,+\\
+\sum^n_{q=1}\sum^n_{r=1}\sum^n_{s=1}\nabla_{\!k}L_q\,g^{rq}\,
L^s\,\nabla_{\!r}L_s-\sum^n_{q=1}\sum^n_{r=1}\nabla_{\!k}L_q\,g^{rq}\,
F_r.
\endgathered
\tag7.24
$$
Further we need to transform $\nabla_{\!k}Q_r$. Using \thetag{3.5},
\thetag{7.2}, and \thetag{7.5}, we get
$$
\hskip -2em
\gathered
\nabla_{\!k}Q_r\compos\lambda=\nabla_{\!k}F_r
+\sum^n_{m=1}v^m\,\nabla_{\!k}\!\nabla_{\!m}L_r\,-\\
-\sum^n_{q=1}\sum^n_{s=1}\nabla_{\!k}L_q\,g^{sq}\,
\tilde\nabla_{\!s}F_r-\sum^n_{q=1}\sum^n_{s=1}
\nabla_{\!k}L_q\,g^{sq}\,\nabla_{\!s}L_r\,-\\
-\sum^n_{q=1}\sum^n_{m=1}\sum^n_{s=1}\nabla_{\!k}L_q\,g^{sq}\,
v^m\,\tilde\nabla_{\!s}\!\nabla_{\!m}L_r.
\endgathered
\tag7.25
$$
Now let's substitute \thetag{7.25} into fourth term in right hand side
of \thetag{5.10}. This yields
$$
\gathered
\sum^n_{r=1}(\nabla_{\!k}Q_r\,W^r)\compos\lambda=\sum^n_{r=1}
L^r\,\nabla_{\!k}F_r+\sum^n_{m=1}\sum^n_{r=1}L^r\,v^m\,\nabla_{\!k}
\!\nabla_{\!m}L_r\,-\\
-\sum^n_{q=1}\sum^n_{r=1}\sum^n_{s=1}L^r\,\nabla_{\!k}L_q\,g^{sq}\,
\tilde\nabla_{\!s}F_r-\sum^n_{q=1}\sum^n_{r=1}\sum^n_{s=1}
L^r\,\nabla_{\!k}L_q\,g^{sq}\,\nabla_{\!s}L_r\,-\\
-\sum^n_{q=1}\sum^n_{m=1}\sum^n_{r=1}\sum^n_{s=1}L^r\,
\nabla_{\!k}L_q\,g^{sq}\,v^m\,\tilde\nabla_{\!s}\!\nabla_{\!m}L_r.
\endgathered\quad
\tag7.26
$$
Next steps are related to curvature tensors in formula \thetag{5.10}.
Curvature tensor $\bold R$ in $\bold v$-representation is given by the
following formula (see Chapter~\uppercase\expandafter{\romannumeral 3}
of thesis \cite{16}):
$$
\hskip -2em
\aligned
R^k_{rij}&=\frac{\partial\Gamma^k_{jr}}{\partial x^i}-
\frac{\partial\Gamma^k_{ir}}{\partial x^j}+\sum^n_{m=1}
\Gamma^k_{im}\,\Gamma^m_{jr}-\sum^n_{m=1}\Gamma^k_{jm}\,
\Gamma^m_{ir}-\\
\vspace{1ex}
&-\sum^n_{m=1}\sum^n_{s=1}v^s\,\Gamma^m_{is}\,\frac{\partial
\Gamma^k_{jr}}{\partial v^m}+\sum^n_{m=1}\sum^n_{s=1}v^s\,
\Gamma^m_{js}\,\frac{\partial\Gamma^k_{ir}}{\partial v^m}.
\endaligned
\tag7.27
$$
Like \thetag{5.11} and \thetag{7.16}, formulas \thetag{5.12}
and \thetag{7.26} express two different extended tensor field
in different representations. However, the relation of these
two tensor fields can be described by the following rather
simple formula:
$$
\pagebreak
R^k_{rij}\compos\lambda=R^k_{rij}+\sum^n_{q=1}\sum^n_{s=1}
\nabla_{\!i}L_q\,g^{sq}\,D^k_{jrs}-\sum^n_{q=1}\sum^n_{s=1}
\nabla_{\!j}L_q\,g^{sq}\,D^k_{irs}.
\quad
\tag7.28
$$
Components of curvature tensor $\bold R$ in right hand side of
\thetag{7.28} are calculated according to the formula \thetag{7.27}.
Using \thetag{7.17} and \thetag{7.28}, we can transform last two
terms in right hand side of \thetag{5.10}. As a result we get
$$
\gather
\hskip -2em
\gathered
-\sum^n_{r=1}\sum^n_{s=1}\sum^n_{m=1}(R^s_{rmk}\,V^m\,W^r\,p_s)
\compos\lambda=-\sum^n_{r=1}\sum^n_{s=1}\sum^n_{m=1}R^s_{rmk}
\,v^m\,L^r\,L_s\,-\\
-\sum^n_{q=1}\sum^n_{r=1}\sum^n_{s=1}\sum^n_{m=1}\sum^n_{a=1}
\nabla_{\!m}L_q\,g^{aq}\,D^s_{kra}\,v^m\,L^r\,L_s\,+\\
+\sum^n_{q=1}\sum^n_{r=1}\sum^n_{s=1}\sum^n_{m=1}\sum^n_{a=1}
\nabla_{\!k}L_q\,g^{aq}\,D^s_{mra}\,v^m\,L^r\,L_s,
\endgathered
\tag7.29
\\
\hskip -2em
\gathered
\sum^n_{r=1}\sum^n_{s=1}\sum^n_{m=1}(D^{sm}_{rk}\,Q_m\,W^r\,p_s)
\compos\lambda=\sum^n_{r=1}\sum^n_{s=1}\sum^n_{m=1}\sum^n_{a=1}
g^{am}\,D^s_{rka}\,\times\\
\times\,F_m\,L^r\,L_s+\sum^n_{q=1}\sum^n_{r=1}\sum^n_{s=1}
\sum^n_{m=1}\sum^n_{a=1}g^{am}\,D^s_{rka}\,v^q\,\nabla_{\!q}L_m
\,L^r\,L_s.
\endgathered
\tag7.30
\endgather
$$
Now we can summarize above calculations and write formula
for components of covector field $\boldsymbol\beta$. Combining
\thetag{7.21}, \thetag{7.23}, \thetag{7.24}, \thetag{7.26},
\thetag{7.29}, and \thetag{7.30}, we get
$$
\gathered
\beta_k=\sum^n_{r=1}v^r\,\nabla_{\!r}U_k+\sum^n_{q=1}F^q\,
\tilde\nabla_{\!q}U_k-\sum^n_{q=1}\sum^n_{r=1}\sum^n_{s=1}
\nabla_{\!k}L_q\,g^{rq}\,v^s\,\nabla_{\!s}L_r\,-\\
-\sum^n_{q=1}\sum^n_{r=1}\nabla_{\!k}L_q\,g^{rq}\,F_r
+\sum^n_{r=1}L^r\,\nabla_{\!k}F_r+\sum^n_{m=1}\sum^n_{r=1}
L^r\,v^m\,\nabla_{\!k}\!\nabla_{\!m}L_r\,-\\
-\sum^n_{q=1}\sum^n_{r=1}\sum^n_{s=1}L^r\,\nabla_{\!k}L_q\,
g^{sq}\,\tilde\nabla_{\!s}F_r-\sum^n_{q=1}\sum^n_{m=1}
\sum^n_{r=1}\sum^n_{s=1}L^r\,\nabla_{\!k}L_q\,g^{sq}\,v^m\,
\times\\
\times\,\tilde\nabla_{\!s}\!\nabla_{\!m}L_r-\sum^n_{r=1}
\sum^n_{s=1}\sum^n_{m=1}R^s_{rmk}\,v^m\,L^r\,L_s+\sum^n_{q=1}
\sum^n_{r=1}\sum^n_{s=1}\sum^n_{m=1}\sum^n_{a=1}g^{aq}\,
\times\\
\times\,\nabla_{\!k}L_q\,D^s_{mra}\,v^m\,L^r\,L_s+\sum^n_{r=1}
\sum^n_{s=1}\sum^n_{m=1}\sum^n_{a=1}g^{am}\,D^s_{rka}\,F_m\,
L^r\,L_s.
\endgathered\quad
\tag7.31
$$
Components of covector field $\bold L$ in \thetag{7.31} are determined
by generalized Legendre transformation \thetag{1.6}, components of its
dual vector field $\bold L$ are given first formula \thetag{7.6}, and
components of covector field $\bold U$ in \thetag{7.31} are determined
by formula \thetag{7.8}. Also we should keep in mind formulas
\thetag{7.7} binding vectorial and covectorial form of force field
$\bold F$ and formulas \thetag{7.16} and \thetag{7.27} for curvature
tensors.\par
    Since $\boldsymbol\alpha$ and $\boldsymbol\beta$ are already
transformed to $\bold v$-representation, formula for components of
covector field $\boldsymbol\eta$ is very simple. It is derived from
\thetag{5.13}:
$$
\pagebreak
\hskip -2em
\eta_k=\beta_k-\sum^n_{s=1}\frac{U_k\,\alpha^s\,L_s}{|\bold L|^2}.
\tag7.32
$$
Here $|\bold L|^2$ is $\bold v$-representation of scalar field
$\Omega$ (see formula \thetag{6.2}). As for weak normality
equations, in $\bold v$-representation we can use their initial
form \thetag{5.14} since all extended tensor fields forming
these equations are already transformed to $\bold v$-re\-presentation
(see \thetag{6.4}, \thetag{7.19}, and \thetag{7.32}).
\head
8. Additional normality equations in $\bold v$-representation.
\endhead
    Transformation of additional normality equations to
$\bold v$-representation consists in transforming extended tensor
fields \thetag{5.15}, \thetag{5.16}, and \thetag{5.17}. Components
of projector field $\bold P$, which enter to \thetag{5.18},
\thetag{5.19}, \thetag{5.20}, and \thetag{5.22}, are already
transformed to $\bold v$-representation (see formulas \thetag{6.3}
or \thetag{6.4}). For $A^{rs}$ we have
$$
\hskip -2em
A^{rs}=\sum^n_{q=1}g^{qr}\,\tilde\nabla_{\!q}L^s.
\tag8.1
$$
In deriving \thetag{8.1} we used \thetag{2.6} and \thetag{6.1}.
Components of tensor field $\bold B$ are given by a little bit more
complicated formula. For them we have
$$
\hskip -2em
\gathered
B^r_s=\sum^n_{q=1}g^{qr}\,\tilde\nabla_{\!q}U_s
+\sum^n_{m=1}\sum^n_{k=1}\sum^n_{q=1}g^{qr}\,L^k
\,L_m\,D^m_{ksq}-\nabla_{\!s}L^r\,+\\
+\sum^n_{q=1}\sum^n_{k=1}\nabla_{\!s}L_k\,g^{qk}\,
\tilde\nabla_{\!q}L^r+\sum^n_{q=1}\sum^n_{m=1}
\frac{g^{qm}\tilde\nabla_{\!q}L^r-g^{qr}\tilde
\nabla_{\!q}L^m}{|\bold L|^2}\,\,U_sL_m.
\endgathered
\tag8.2
$$
In deriving \thetag{8.2} we used formulas \thetag{7.22}, \thetag{7.17},
\thetag{7.10}, and \thetag{2.6}. In the last step we should transform
formula \thetag{5.17} for components of extended tensor field $\bold C$
to $\bold v$-representation. Using \thetag{7.20}, \thetag{7.22},
\thetag{7.10}, \thetag{7.17}, and \thetag{7.28}, we get
$$
\hskip -2em
\gathered
C_{rs}=\nabla_{\!r}U_s-\sum^n_{q=1}\sum^n_{k=1}\nabla_{\!r}L_q\,
g^{kq}\,\tilde\nabla_{\!k}U_s
-\sum^n_{m=1}\frac{U_s\,\nabla_{\!r}L^m}{|\bold L|^2}\,L_m\,-\\
-\sum^n_{q=1}\sum^n_{m=1}\frac{U_r\,g^{qm}\,\tilde\nabla_{\!q}U_s}
{|\bold L|^2}\,L_m
+\sum^n_{m=1}\sum^n_{q=1}\sum^n_{k=1}\frac{U_s\,\nabla_{\!r}L_k
\,g^{qk}\,\tilde\nabla_{\!q}L^m}{|\bold L|^2}\,L_m\,-\\
-\sum^n_{k=1}\sum^n_{q=1}\sum^n_{m=1}\sum^n_{a=1}
\frac{g^{aq}\,D^m_{ksa}\,U_r}{|\bold L|^2}\,L_m\,L^kL_q
-\sum^n_{k=1}\sum^n_{q=1}\frac{R^q_{krs}}{2}\,L^kL_q\,-\\\
-\sum^n_{k=1}\sum^n_{q=1}\sum^n_{m=1}\sum^n_{a=1}
\nabla_{\!r}L_m\,D^q_{ksa}\,g^{am}\,L^kL_q.
\endgathered
\tag8.3
$$
Substituting \thetag{8.1}, \thetag{8.2}, \thetag{8.3} into
\thetag{5.18}, \thetag{5.19}, \thetag{5.20}, we get additional
normality equations transformed to $\bold v$-representation.
Formula \thetag{5.21} for scalar parameter $\lambda$ remains
unchanged in $\bold v$-representation as well.
\head
9. Connection invariance.
\endhead
    Let's recall that \thetag{4.5} is a tensorial form of
differential equations \thetag{1.1}. In order to write
\thetag{4.5} we need to choose some extended connection
$\Gamma$ in $M$. However, strong normality condition in
definition~5.2, from which weak and additional normality
equation were derived in \cite{26}, is irrelative to the
choice of $\Gamma$ when applied to Newtonian dynamical
system \thetag{1.1}. Therefore we should prove that
complete system of normality equations is invariant under
the change of connection components, provided functions
$\Phi^i$ in differential equations \thetag{1.1} are
unchanged. This means that we should consider gauge
transformations
$$
\align
&\hskip -2em
\Gamma^k_{ij}\to\Gamma^k_{ij}+T^k_{ij},\\
\vspace{-1ex}
&\hskip -2em
\tag9.1\\
\vspace{-1ex}
&\hskip -2em
F^i\to F^i+\sum^n_{k=1}\sum^n_{j=1}T^i_{kj}\,v^k\,v^j.
\endalign
$$
Here $T^k_{ij}$ are components of symmetric extended tensor field
$\bold T$. Tensor fields $\bold L$, $\bold g$, $\Omega=|\bold L|^2$,
and $\bold P$ are obviously invariant under gauge transformations
\thetag{9.1}:
$$
\xalignat 2
&\hskip -2em
L_i\to L_i, &&g_{ij}\to g_{ij},\\
&\hskip -2em
g^{ij}\to g^{ij}, &&L^i\to L^i,
\tag9.2\\
&\hskip -2em
|\bold L|^2\to |\bold L|^2, &&P^i_k\to P^i_k.
\endxalignat
$$
Applying \thetag{9.2} to \thetag{8.1}, we find that extended tensor
field $\bold A$ with components \thetag{8.1} is invariant under
gauge transformations \thetag{9.1}:
$$
\hskip -2em
A^{rs}\to A^{rs}.
\tag9.3
$$
Due to \thetag{9.3} and \thetag{9.2} additional normality equation
\thetag{5.18} is invariant under gauge transformations \thetag{9.1}.
\par
     Applying formulas \thetag{9.1} and \thetag{9.2} to \thetag{7.7},
for covectorial form of force field of dynamical system \thetag{4.5}
we derive the following transformation rule:
$$
\hskip -2em
F_i\to F_i+\sum^n_{s=1}\sum^n_{k=1}\sum^n_{j=1}g_{is}
\,T^s_{kj}\,v^k\,v^j.
\tag9.4
$$
Connection components are used in formula \thetag{3.4}. Therefore
we have
$$
\hskip -2em
\aligned
&\nabla_{\!s}L_k\to \nabla_{\!s}L_k-\sum^n_{a=1}T^a_{sk}\,L_a
-\sum^n_{a=1}\sum^n_{b=1}g_{kb}\,T^b_{as}\,v^a,\\
&\nabla_{\!s}L^k\to \nabla_{\!s}L^k+\sum^n_{a=1}T^k_{sa}\,L^a
-\sum^n_{a=1}\sum^n_{b=1}v^a\,T^b_{as}\,\tilde\nabla_{\!b}L^k.
\endaligned
\tag9.5
$$
From \thetag{9.4} and \thetag{9.5} for covector field $\bold U$
we derive transformation rule
$$
\hskip -2em
U_i\to U_i+\sum^n_{q=1}\sum^n_{r=1}L^q\,T^r_{iq}\,L_r.
\tag9.6
$$
For dynamic curvature tensor $\bold D$ by differentiating \thetag{9.1}
we get
$$
\hskip -2em
D^k_{rij}\to D^k_{rij}-\tilde\nabla_{\!j}T^k_{ir}.
\tag9.7
$$
Curvature tensor $\bold R$ is determined by connection components
according to formula \thetag{7.28}. Here we have the following
transformation rule:
$$
\hskip -2em
\aligned
&R^k_{rij}\to R^k_{rij}+\nabla_{\!i}T^k_{jr}-\nabla_{\!j}T^k_{ir}
-\sum^n_{s=1}\sum^n_{m=1}v^m\,T^s_{jm}\,D^k_{irs}\,+\\
&+\,\sum^n_{s=1}\sum^n_{m=1}v^m\,T^s_{im}\,D^k_{jrs}
+\sum^n_{m=1}\left(T^k_{im}\,T^m_{jr}-T^k_{jm}\,T^m_{ir}\right)+\\
&+\,\sum^n_{s=1}\sum^n_{m=1}v^m\,T^s_{jm}\,\tilde\nabla_{\!s}T^k_{ir}
-\sum^n_{s=1}\sum^n_{m=1}v^m\,T^s_{im}\,\tilde\nabla_{\!s}T^k_{jr}.
\endaligned
\tag9.8
$$
Further let's derive transformation rules for covariant derivatives
of covector field $\bold U$. For covariant derivative $\tilde
\nabla_{\!q}U_s$, applying formula \thetag{9.6}, we obtain
$$
\hskip -2em
\gathered
\tilde\nabla_{\!q}U_s\to\tilde\nabla_{\!q}U_s+
\sum^n_{k=1}\sum^n_{r=1}\tilde\nabla_{\!q}L^k\,T^r_{sk}
\,L_r\,+\\
+\sum^n_{k=1}\sum^n_{r=1}L^k\,\tilde\nabla_{\!q}T^r_{sk}\,L_r+
\sum^n_{k=1}\sum^n_{r=1}L^k\,T^r_{sk}\,\tilde\nabla_{\!q}L_r.
\endgathered
\tag9.9
$$
Transformation rule for covariant derivative $\nabla_{\!q}U_s$ is
more complicated:
$$
\gathered
\nabla_{\!q}U_s\to\nabla_{\!q}U_s+\sum^n_{k=1}\sum^n_{r=1}
\nabla_{\!q}L^k\,T^r_{sk}\,L_r+\sum^n_{k=1}\sum^n_{r=1}
L^k\,\nabla_{\!q}T^r_{sk}\,L_r\,+\\
+\sum^n_{k=1}\sum^n_{r=1}L^k\,T^r_{sk}\,\nabla_{\!q}L_r
-\sum^n_{m=1}T^m_{qs}\,U_m-\sum^n_{m=1}\sum^n_{k=1}\sum^n_{r=1}
T^m_{qs}\,L^k\,T^r_{mk}\,L_r\,-\\
-\sum^n_{m=1}\sum^n_{a=1}v^a\,T^m_{aq}\,\tilde\nabla_{\!m}U_s
-\sum^n_{m=1}\sum^n_{a=1}\sum^n_{k=1}\sum^n_{r=1}v^a\,T^m_{aq}\,
\tilde\nabla_{\!m}L^k\,T^r_{sk}\,L_r\,-\\
-\sum^n_{m=1}\sum^n_{a=1}v^a\,T^m_{aq}\left(\,\shave{\sum^n_{k=1}
\sum^n_{r=1}}L^k\,\tilde\nabla_{\!m}T^r_{sk}\,L_r
+\shave{\sum^n_{k=1}\sum^n_{r=1}}L^k\,T^r_{sk}\,\tilde\nabla_{\!m}
L_r\!\right).
\endgathered\quad
\tag9.10
$$
Now we can combine \thetag{9.9}, \thetag{9.7}, \thetag{9.6}, and
\thetag{9.5} according to formula \thetag{8.2}:
$$
\gather
B^r_s\to B^r_s+\sum^n_{q=1}\sum^n_{k=1}\sum^n_{m=1}g^{qr}\,\tilde
\nabla_{\!q}L^k\,T^m_{sk}\,L_m+\sum^n_{q=1}\sum^n_{k=1}\sum^n_{m=1}
g^{qr}\,L^k\,\tilde\nabla_{\!q}T^m_{sk}\,L_m\,+\\
+\sum^n_{q=1}\sum^n_{k=1}\sum^n_{m=1}g^{qr}\,L^k\,T^m_{sk}\,\tilde
\nabla_{\!q}L_m-\sum^n_{m=1}\sum^n_{k=1}\sum^n_{q=1}g^{qr}\,L^k
\,L_m\,\tilde\nabla_{\!q}T^m_{sk}-\sum^n_{m=1}T^r_{sm}\,L^m\,+\\
+\sum^n_{a=1}\sum^n_{b=1}v^a\,T^b_{as}\,\tilde\nabla_{\!b}L^r
-\sum^n_{q=1}\sum^n_{k=1}\sum^n_{m=1}T^m_{sk}\,L_m\,g^{qk}\,
\tilde\nabla_{\!q}L^r-\sum^n_{q=1}\sum^n_{k=1}\sum^n_{a=1}
\sum^n_{b=1}g_{kb}\,\times\\
\times\,T^b_{as}\,v^a\,g^{qk}\,\tilde\nabla_{\!q}L^r+\sum^n_{q=1}
\sum^n_{k=1}\sum^n_{m=1}\frac{A^{mr}-A^{rm}}{|\bold L|^2}\,L_m
\,L^q\,T^k_{sq}\,L_k.
\endgather
$$
Using \thetag{2.7} for $\tilde\nabla_{\!q}L_m$ and canceling
similar terms in the above formula, we get
$$
\hskip -2em
B^r_s\to B^r_s+\sum^n_{m=1}\sum^n_{q=1}\sum^n_{k=1}L_m\,T^m_{sq}
\,P^q_k\left(A^{rk}-A^{kr}\right).
\tag9.11
$$
Here we used formula \thetag{6.4} for $P^q_k$. Substituting
\thetag{9.11} into normality equation \thetag{5.19} and taking
into account another normality equation \thetag{5.18}, we find
that \thetag{5.19} is invariant under gauge transformations
\thetag{9.1}.\par
    Now let's apply transformation rules \thetag{9.10}, \thetag{9.9},
\thetag{9.8}, \thetag{9.7}, \thetag{9.6}, and \thetag{9.5} to formula
\thetag{8.3} for components of covector field $\bold C$ in
$\bold v$-representation:
$$
\allowdisplaybreaks
\gather
C_{rs}\to C_{rs}+\ldots+\sum^n_{k=1}\sum^n_{q=1}
\nabla_{\!r}L^k\,T^q_{sk}\,L_q+\sum^n_{k=1}\sum^n_{q=1}
L^k\,\nabla_{\!r}T^q_{sk}\,L_q+\sum^n_{k=1}\sum^n_{q=1}
L^k\,\times\\
\times\,T^q_{sk}\,\nabla_{\!r}L_q
-\sum^n_{m=1}\sum^n_{a=1}v^a\,T^m_{ar}\,\tilde\nabla_{\!m}U_s
-\sum^n_{m=1}\sum^n_{a=1}v^a\,T^m_{ar}\sum^n_{k=1}\sum^n_{q=1}
\left(\tilde\nabla_{\!m}L^k\,T^q_{sk}\,L_q\,+\right.\\
\left.+\,L^k\,\tilde\nabla_{\!m}T^q_{sk}\,L_q+L^k\,T^q_{sk}\,
\tilde\nabla_{\!m}L_q\!\right)-\sum^n_{k=1}\sum^n_{q=1}
\sum^n_{b=1}\sum^n_{a=1}\nabla_{\!r}L_q\,g^{kq}\,
\tilde\nabla_{\!k}L^b\,T^a_{sb}\,L_a\,-\\
-\sum^n_{k=1}\sum^n_{q=1}\sum^n_{b=1}\sum^n_{a=1}\nabla_{\!r}L_q\,g^{kq}
\,L^b\,\tilde\nabla_{\!k}T^a_{sb}\,L_a
-\sum^n_{k=1}\sum^n_{q=1}\sum^n_{b=1}\sum^n_{a=1}\nabla_{\!r}L_q\,g^{kq}
\,L^b\,T^a_{sb}\,\times\\
\times\,\tilde\nabla_{\!k}L_a+\sum^n_{k=1}\sum^n_{q=1}
\sum^n_{m=1}T^m_{rq}\,L_m\,g^{kq}\,\tilde\nabla_{\!k}U_s
+\sum^n_{k=1}\sum^n_{q=1}\sum^n_{m=1}\sum^n_{b=1}\sum^n_{a=1}
T^m_{rq}\,L_m\,g^{kq}\,\times\\
\times\,\tilde\nabla_{\!k}L^b\,T^a_{sb}\,L_a+\sum^n_{b=1}
\sum^n_{a=1}\sum^n_{k=1}\sum^n_{q=1}\sum^n_{m=1}
T^m_{rq}\,L_m\,g^{kq}\,L^b\,\tilde\nabla_{\!k}T^a_{sb}\,L_a
+\sum^n_{q=1}\sum^n_{m=1}T^m_{rq}\,\times\\
\times\,L_m\,\sum^n_{b=1}\sum^n_{a=1}\sum^n_{k=1}
g^{kq}\,L^b\,T^a_{sb}\,\tilde\nabla_{\!k}L_a
+\sum^n_{k=1}\sum^n_{m=1}T^k_{mr}\,v^m\,\tilde\nabla_{\!k}U_s
+\sum^n_{k=1}\sum^n_{m=1}T^k_{mr}\,\times\\
\times\,\sum^n_{b=1}\sum^n_{a=1}v^m\,\tilde\nabla_{\!k}L^b
\,T^a_{sb}\,L_a+\sum^n_{k=1}\sum^n_{b=1}\sum^n_{a=1}\sum^n_{m=1}
T^k_{mr}\,v^m\,L^b\,\tilde\nabla_{\!k}T^a_{sb}\,L_a+\sum^n_{m=1}
v^m\,\times\\
\times\,\sum^n_{k=1}\sum^n_{b=1}\sum^n_{a=1}T^k_{mr}\,L^b\,T^a_{sb}
\,\tilde\nabla_{\!k}L_a+\sum^n_{k=1}\sum^n_{q=1}\sum^n_{m=1}
\sum^n_{a=1}\sum^n_{b=1}T^b_{rm}\,L_b\,g^{am}\,D^q_{ksa}\,L^kL_q\,+\\
+\sum^n_{k=1}\sum^n_{q=1}\sum^n_{a=1}\sum^n_{b=1}T^a_{br}\,v^b
\,D^q_{ksa}\,L^kL_q+\sum^n_{k=1}\sum^n_{q=1}\sum^n_{m=1}\sum^n_{a=1}
\nabla_{\!r}L_m\,g^{am}\,\tilde\nabla_{\!a}T^q_{sk}\,L^kL_q\,-\\
-\sum^n_{k=1}\sum^n_{q=1}\sum^n_{m=1}\sum^n_{a=1}\sum^n_{b=1}
T^b_{rm}\,L_b\,g^{am}\,\tilde\nabla_{\!a}T^q_{sk}\,L^kL_q
-\sum^n_{k=1}\sum^n_{q=1}\sum^n_{a=1}\sum^n_{b=1}T^a_{br}\,v^b\,
\tilde\nabla_{\!a}T^q_{sk}\,\times\\
\times\,L^kL_q-\sum^n_{k=1}\sum^n_{q=1}\nabla_{\!r}T^q_{sk}
\,L^kL_q-\sum^n_{k=1}\sum^n_{q=1}\sum^n_{a=1}\sum^n_{m=1}v^m
\,T^a_{rm}\,D^q_{ska}\,L^kL_q-\sum^n_{k=1}L^k\,\times\\
\times\,\sum^n_{q=1}\sum^n_{m=1}T^q_{rm}\,T^m_{sk}\,L_q
+\sum^n_{k=1}\sum^n_{q=1}\sum^n_{a=1}\sum^n_{m=1}v^m\,T^a_{rm}
\,\tilde\nabla_{\!a}T^q_{sk}\,L^kL_q+\frac{1}{|\bold L|^2}
\,Z_{rs}.
\endgather
$$
We used $Z_{rs}$ in order to denote contribution of those terms
in \thetag{8.3} which have denominator $|\bold L|^2$. By dots
we denote terms symmetric in indices $r$ and $s$. They do not affect
the ultimate form of the equations \thetag{5.20}. When collecting
similar terms in the above formula some terms cancel each other.
As a result we have 
$$
\gathered
C_{rs}\to C_{rs}+\ldots+\sum^n_{k=1}\sum^n_{q=1}\nabla_{\!r}L^k
\,T^q_{sk}\,L_q-\sum^n_{k=1}\sum^n_{q=1}\sum^n_{b=1}\sum^n_{a=1}
\nabla_{\!r}L_q\,\times\\
\times\,g^{kq}\,\tilde\nabla_{\!k}L^b\,T^a_{sb}\,L_a+\sum^n_{k=1}
\sum^n_{q=1}\sum^n_{m=1}T^m_{rq}\,L_m\,g^{kq}\,\tilde\nabla_{\!k}U_s
+\sum^n_{k=1}\sum^n_{q=1}g^{kq}\,\times\\
\times\,\sum^n_{m=1}\sum^n_{b=1}\sum^n_{a=1}T^m_{rq}\,L_m
\,\tilde\nabla_{\!k}L^b\,T^a_{sb}\,L_a+\sum^n_{k=1}\sum^n_{q=1}
\sum^n_{m=1}\sum^n_{a=1}\sum^n_{b=1}T^b_{rm}\,\times\\
\vspace{1ex}
\times\,L_b\,g^{am}\,D^q_{ksa}\,L^kL_q+\frac{1}{|\bold L|^2}
\,Z_{rs}.
\endgathered
\tag9.12
$$
In deriving \thetag{9.12} we used formula $g_{qk}=\tilde\nabla_{\!k}
L_q$ (see \thetag{2.7} above). For $Z_{rs}$ in \thetag{9.12} we have
rather huge formula derived from formula \thetag{8.3} for $C_{rs}$:
$$
\allowdisplaybreaks
\gather
Z_{rs}=-\sum^n_{m=1}\sum^n_{a=1}T^m_{ra}\,L^a\,L_m\,U_s+\sum^n_{m=1}
\sum^n_{a=1}\sum^n_{b=1}v^a\,T^b_{ar}\,\tilde\nabla_{\!b}L^m\,L_m
\,U_s\,-\\
-\sum^n_{q=1}\sum^n_{k=1}\sum^n_{m=1}\nabla_{\!r}L^m\,L_m\,L^q
\,T^k_{sq}\,L_k-\sum^n_{q=1}\sum^n_{k=1}\sum^n_{m=1}\sum^n_{a=1}
T^m_{ra}\,L^a\,L_m\,L^q\,T^k_{sq}\,L_k\,+\\
+\sum^n_{q=1}\sum^n_{k=1}\sum^n_{m=1}\sum^n_{a=1}\sum^n_{b=1}
v^a\,T^b_{ar}\,\tilde\nabla_{\!b}L^m\,L_m\,L^q\,T^k_{sq}\,L_k
-\sum^n_{k=1}\sum^n_{q=1}\sum^n_{m=1}\sum^n_{c=1}U_r\,g^{qm}
\,\times\\
\times\,\tilde\nabla_{\!q}L^k\,T^c_{sk}\,L_c\,L_m-\sum^n_{k=1}
\sum^n_{q=1}\sum^n_{m=1}\sum^n_{c=1}U_r\,g^{qm}\,L^k
\,\tilde\nabla_{\!q}T^c_{sk}\,L_c\,L_m-\sum^n_{k=1}L^k\,\times\\
\times\,\sum^n_{q=1}\sum^n_{m=1}\sum^n_{c=1}U_r\,g^{qm}\,T^c_{sk}
\,\tilde\nabla_{\!q}L_c\,L_m-\sum^n_{q=1}\sum^n_{m=1}\sum^n_{a=1}
\sum^n_{b=1}L^a\,T^b_{ra}\,L_b\,g^{qm}\,\tilde\nabla_{\!q}U_s
\,\times\\
\times\,L_m-\sum^n_{k=1}\sum^n_{q=1}\sum^n_{m=1}\sum^n_{a=1}
\sum^n_{b=1}\sum^n_{c=1}L^a\,T^b_{ra}\,L_b\,g^{qm}\,
\tilde\nabla_{\!q}L^k\,T^c_{sk}\,L_c\,L_m-\sum^n_{k=1}L^k
\,\times\\
\times\,\sum^n_{q=1}\sum^n_{m=1}\sum^n_{a=1}\sum^n_{b=1}\sum^n_{c=1}
L^a\,T^b_{ra}\,L_b\,g^{qm}\,\tilde\nabla_{\!q}T^c_{sk}\,L_c\,L_m
-\sum^n_{k=1}\sum^n_{q=1}\sum^n_{m=1}\sum^n_{a=1}L^a\,\times\\
\times\sum^n_{b=1}\sum^n_{c=1}T^b_{ra}\,L_b\,g^{qm}\,L^k\,T^c_{sk}
\,\tilde\nabla_{\!q}L_c\,L_m-\sum^n_{k=1}\sum^n_{q=1}\sum^n_{m=1}
\sum^n_{c=1}U_s\,T^c_{rk}\,L_c\,g^{qk}\,\times\\
\times\,\tilde\nabla_{\!q}L^m\,L_m-\sum^n_{k=1}\sum^n_{q=1}
\sum^n_{m=1}\sum^n_{c=1}\sum^n_{e=1}U_s\,g_{ke}\,T^e_{cr}
\,v^c\,g^{qk}\,\tilde\nabla_{\!q}L^m\,L_m+\sum^n_{a=1}L^a
\,\times\\
\times\,\sum^n_{k=1}\sum^n_{q=1}\sum^n_{m=1}\sum^n_{b=1}T^b_{sa}
\,L_b\,\nabla_{\!r}L_k\,g^{qk}\,\tilde\nabla_{\!q}L^m\,L_m
-\sum^n_{k=1}\sum^n_{q=1}\sum^n_{m=1}\sum^n_{a=1}\sum^n_{b=1}
\sum^n_{c=1}L^a\,\times\\
\times\,T^b_{sa}\,L_b\,T^c_{rk}\,L_c\,g^{qk}\,\tilde\nabla_{\!q}L^m
\,L_m-\sum^n_{k=1}\sum^n_{q=1}\sum^n_{m=1}\sum^n_{a=1}\sum^n_{b=1}
\sum^n_{c=1}\sum^n_{e=1}L^a\,T^b_{sa}\,L_b\,g_{ke}\,\times\\
\times\,T^e_{cr}\,v^c\,g^{qk}\,\tilde\nabla_{\!q}L^m\,L_m+\sum^n_{k=1}
\sum^n_{q=1}\sum^n_{m=1}\sum^n_{a=1}g^{aq}\,\tilde\nabla_{\!a}T^m_{sk}
\,U_r\,L_m\,L^kL_q-\sum^n_{k=1}L^k\,\times\\
\times\,\sum^n_{q=1}\sum^n_{m=1}\sum^n_{a=1}\sum^n_{b=1}\sum^n_{c=1}
g^{aq}\,D^m_{ksa}\,L^b\,T^c_{rb}\,L_c\,L_m\,L_q+\sum^n_{k=1}
\sum^n_{q=1}\sum^n_{m=1}\sum^n_{a=1}\sum^n_{b=1}L^k\times\\
\times\,\sum^n_{c=1}g^{aq}\,\tilde\nabla_{\!a}T^m_{sk}\,L^b
\,T^c_{rb}\,L_c\,L_m\,L_q.
\endgather
$$
Upon collecting and canceling similar terms for $Z_{rs}$ we get
shorter expression
$$
\gathered
Z_{rs}=\dots-\sum^n_{q=1}\sum^n_{k=1}\sum^n_{m=1}\nabla_{\!r}L^m
\,L_m\,L^q\,T^k_{sq}\,L_k-\sum^n_{k=1}\sum^n_{q=1}
\sum^n_{m=1}U_r\,\times\\
\times\,\sum^n_{c=1}g^{qm}\,\tilde\nabla_{\!q}L^k\,T^c_{sk}
\,L_c\,L_m-\sum^n_{q=1}\sum^n_{m=1}\sum^n_{a=1}\sum^n_{b=1}
L^a\,T^b_{ra}\,L_b\,g^{qm}\,\times\\
\times\,\tilde\nabla_{\!q}U_s\,L_m-\sum^n_{a=1}\sum^n_{b=1}
\sum^n_{k=1}\sum^n_{q=1}\sum^n_{m=1}\sum^n_{c=1}T^b_{ra}\,
L^a\,L_b\,g^{qm}\,\tilde\nabla_{\!q}L^k\,T^c_{sk}\,\times\\
\times\,L_c\,L_m-\sum^n_{k=1}\sum^n_{q=1}\sum^n_{m=1}
\sum^n_{c=1}U_s\,T^c_{rk}\,L_c\,g^{qk}\,\tilde\nabla_{\!q}L^m
\,L_m+\sum^n_{a=1}L^a\,\times\\
\times\,\sum^n_{k=1}\sum^n_{q=1}\sum^n_{m=1}\sum^n_{b=1}T^b_{sa}
\,L_b\,\nabla_{\!r}L_k\,g^{qk}\,\tilde\nabla_{\!q}L^m\,L_m
-\sum^n_{a=1}\sum^n_{k=1}\sum^n_{q=1}L^a\,\times\\
\times\,\sum^n_{m=1}\sum^n_{b=1}\sum^n_{c=1}T^b_{sa}\,L_b
\,T^c_{rk}\,L_c\,g^{qk}\,\tilde\nabla_{\!q}L^m\,L_m-\sum^n_{a=1}
\sum^n_{k=1}\sum^n_{q=1}g^{aq}\,\times\\
\times\sum^n_{m=1}\sum^n_{b=1}\sum^n_{c=1}D^m_{ksa}\,L^b\,T^c_{rb}
\,L_c\,L_m\,L^kL_q.
\endgathered\quad
\tag9.13
$$
Now we can substitute \thetag{9.13} into \thetag{9.12} for further
transformations. Keeping in mind formula \thetag{6.4} for components
of projector field $\bold P$, we get
$$
\allowdisplaybreaks
\gather
C_{rs}\to C_{rs}+\ldots+\sum^n_{q=1}\sum^n_{b=1}\sum^n_{m=1}
T^b_{rq}\,L_b\,P^q_m\left(\,\shave{\sum^n_{k=1}}g^{km}\,\tilde
\nabla_{\!k}U_s+\shave{\sum^n_{k=1}\sum^n_{a=1}\sum^n_{c=1}}
g^{am}\,\times\right.\\
\left.\times\,D^c_{ksa}\,L^k\,L_c-\nabla_{\!s}L^m+\sum^n_{k=1}
\sum^n_{c=1}g^{kc}\,\tilde\nabla_{\!k}L^m\,\nabla_{\!s}L_c\!\right)
+\sum^n_{q=1}\sum^n_{b=1}\sum^n_{m=1}T^b_{rq}\,L_b\,L_m\,U_s
\times\\
\times\,\sum^n_{k=1}\frac{g^{km}\,\tilde\nabla_{\!k}L^q-
g^{kq}\,\tilde\nabla_{\!k}L^m}{|\bold L|^2}+\sum^n_{q=1}\sum^n_{b=1}
\sum^n_{m=1}\sum^n_{k=1}\sum^n_{a=1}\sum^n_{c=1}T^b_{rq}\,L_b\,P^q_m
\,g^{km}\,\tilde\nabla_{\!k}L^a\,\times\\
\times\,T^c_{sa}\,L_c-\sum^n_{q=1}\sum^n_{b=1}\sum^n_{m=1}
\sum^n_{k=1}\sum^n_{a=1}\sum^n_{c=1}T^b_{rq}\,L_b\,
\frac{g^{kq}\,\tilde\nabla_{\!k}L^m\,L_m\,T^a_{sc}\,L_a\,L^c}
{|\bold L|^2}.
\endgather
$$
One can find four summands in right hand side of this relationship
for $C_{rs}$ (sums in round brackets are considered as a single
term). Let's transform second summand containing fraction with
numerator skew-symmetric in indices $q$ and $ m$:
$$
\gather
\sum^n_{m=1}L_m\,\frac{g^{km}\,\tilde\nabla_{\!k}L^q-
g^{kq}\,\tilde\nabla_{\!k}L^m}{|\bold L|^2}=
\sum^n_{m=1}\sum^n_{a=1}P^q_a\,L_m\,\frac{g^{km}
\,\tilde\nabla_{\!k}L^a-g^{ka}\,\tilde\nabla_{\!k}L^m}
{|\bold L|^2}\,+\\
+\sum^n_{m=1}\sum^n_{a=1}\frac{L^q\,L_a}{|\bold L|^2}
\,L_m\,\frac{g^{km}\,\tilde\nabla_{\!k}L^a-g^{ka}\,
\tilde\nabla_{\!k}L^m}{|\bold L|^2}.
\endgather
$$
Here we used formula \thetag{6.4} for $P^q_a$ again. Note that
last term in right hand side of the above formula vanishes due
to skew symmetry of the expression under summation with respect
to indices $m$ and $a$. Then for $C_{rs}$ we obtain
$$
\allowdisplaybreaks
\gather
C_{rs}\to C_{rs}+\ldots+\sum^n_{q=1}\sum^n_{b=1}\sum^n_{m=1}
T^b_{rq}\,L_b\,P^q_m\,B^m_s+\sum^n_{q=1}\sum^n_{b=1}
\sum^n_{m=1}\sum^n_{k=1}\sum^n_{a=1}\sum^n_{c=1}T^b_{rq}\,L_b
\,\times\\
\times\,P^q_m\,g^{km}\,\tilde\nabla_{\!k}L^a\,T^c_{sa}\,L_c
-\sum^n_{q=1}\sum^n_{b=1}\sum^n_{m=1}\sum^n_{k=1}\sum^n_{a=1}
\sum^n_{c=1}T^b_{rq}\,L_b\,\frac{g^{kq}\,\tilde\nabla_{\!k}L^m
\,L_m\,T^a_{sc}\,L_a\,L^c}{|\bold L|^2}.
\endgather
$$
In deriving this formula we used \thetag{8.2} for $B^m_s$. Now we
shall apply similar trick with skew symmetry for transforming last
term in the above formula:
$$
\gathered
\sum^n_{q=1}\sum^n_{b=1}\sum^n_{m=1}\sum^n_{k=1}\sum^n_{a=1}
\sum^n_{c=1}T^b_{rq}\,L_b\,\frac{g^{kq}\,\tilde\nabla_{\!k}L^m
\,L_m\,T^a_{sc}\,L_a\,L^c}{|\bold L|^2}=\\
\sum^n_{q=1}\sum^n_{b=1}\sum^n_{m=1}\sum^n_{k=1}\sum^n_{a=1}
\sum^n_{c=1}\sum^n_{e=1}T^b_{rq}\,L_b\,P^q_e\,\frac{g^{ke}\,
\tilde\nabla_{\!k}L^m\,L_m\,T^a_{sc}\,L_a\,L^c}{|\bold L|^2}\,+\\
+\sum^n_{q=1}\sum^n_{b=1}\sum^n_{m=1}\sum^n_{k=1}\sum^n_{a=1}
\sum^n_{c=1}\sum^n_{e=1}T^b_{rq}\,L_b\,\frac{L^q\,L_e}
{|\bold L|^2}\,\frac{g^{ke}\,\tilde\nabla_{\!k}L^m\,L_m
\,T^a_{sc}\,L_a\,L^c}{|\bold L|^2}.
\endgathered\quad
\tag9.14
$$
Last term in \thetag{9.14} is symmetric with respect to indices
$r$ and $s$. Therefore we can denote it by dots when substituting
\thetag{9.14} into formula for $C_{rs}$. As a result we get
$$
\hskip -2em
\gathered
C_{rs}\to C_{rs}+\ldots+\sum^n_{q=1}\sum^n_{b=1}\sum^n_{m=1}
T^b_{rq}\,L_b\,P^q_m\,B^m_s\,+\\
+\sum^n_{q=1}\sum^n_{b=1}\sum^n_{m=1}\sum^n_{a=1}\sum^n_{c=1}
\sum^n_{e=1}T^b_{rq}\,L_b\,P^q_m\,A^{ma}\,P^c_a\,T^e_{sc}\,L_e.
\endgathered
\tag9.15
$$
Here $A^{ma}$ are given by formula \thetag{8.1}. Substituting
\thetag{9.15} into the equation \thetag{5.20} and taking into
account equations \thetag{5.18} and \thetag{5.19}, we get
$$
\pagebreak
\sum^n_{r=1}\sum^n_{s=1}(C_{rs}-C_{sr})\,P^r_i\,P^s_j
\to\sum^n_{r=1}\sum^n_{s=1}(C_{rs}-C_{sr})\,P^r_i\,P^s_j.
$$
This means that the equation \thetag{5.20} is invariant under
gauge transformations \thetag{9.1}. Summarizing this result and
similar results for \thetag{5.18} and \thetag{5.19} obtained
above, we can formulate them in the following theorem.
\proclaim{Theorem 9.1} Additional normality equations \thetag{5.18},
\thetag{5.19}, and \thetag{5.20} transformed to
$\bold v$-representation are invariant under gauge transformations
\thetag{9.1}.
\endproclaim
    Now let's consider weak normality equations \thetag{5.14}.
In order to prove similar theorem for them we should derive
transformation rules for vector field $\boldsymbol\alpha$ and
covector field $\boldsymbol\eta$. For applying \thetag{9.1}
to \thetag{7.19} we need to perform some preliminary calculations.
From \thetag{9.5} for covariant derivative $\tilde\nabla_{\!q}
\nabla_{\!s}L_r$ we derive
$$
\gathered
\tilde\nabla_{\!q}\nabla_{\!s}L_r\to\tilde\nabla_{\!q}\nabla_{\!s}L_r
-\sum^n_{a=1}\tilde\nabla_{\!q}T^a_{sr}\,L_a-\sum^n_{a=1}T^a_{sr}
\,g_{aq}\,-\\
-\sum^n_{a=1}\sum^n_{b=1}\tilde\nabla_{\!q}g_{rb}\,T^b_{as}\,v^a
-\sum^n_{a=1}\sum^n_{b=1}g_{rb}\,\tilde\nabla_{\!q}T^b_{as}\,v^a
-\sum^n_{b=1}g_{rb}\,T^b_{qs}.
\endgathered
\tag9.16
$$
Now, using \thetag{7.19}, \thetag{9.4}, \thetag{9.5}, \thetag{9.6},
\thetag{9.7}, and \thetag{9.16}, for $\alpha^k$ we obtain
$$
\gather
\alpha^k\to\alpha^k+\sum^n_{r=1}\sum^n_{q=1}\sum^n_{a=1}g^{rk}\,L^q
\,T^a_{rq}\,L_a+\sum^n_{r=1}\sum^n_{a=1}v^r\,T^k_{ra}\,L^a
-\sum^n_{r=1}\sum^n_{a=1}\sum^n_{b=1}v^r\,v^a\,\times\\
\times\,T^b_{ar}\,
\tilde\nabla_{\!b}L^k+\sum^n_{q=1}\sum^n_{a=1}\sum^n_{j=1}T^q_{aj}
\,v^a\,v^j\,\tilde\nabla_qL^k-\sum^n_{q=1}\sum^n_{r=1}\sum^n_{s=1}
\sum^n_{a=1}L^r\,g^{qk}\,v^s\,\tilde\nabla_{\!q}T^a_{sr}\,L_a\,-\\
-\sum^n_{r=1}\sum^n_{s=1}
L^r\,v^s\,T^k_{sr}-\sum^n_{q=1}\sum^n_{r=1}\sum^n_{s=1}\sum^n_{a=1}
\sum^n_{b=1}L^r\,g^{qk}\,v^s\,\tilde\nabla_{\!q}g_{rb}\,T^b_{as}
\,v^a-\sum^n_{a=1}\sum^n_{b=1}L_b\,\times\\
\times\,\sum^n_{q=1}\sum^n_{s=1}\,g^{qk}\,v^s
\,\tilde\nabla_{\!q}T^b_{as}\,v^a-\sum^n_{q=1}\sum^n_{s=1}
\sum^n_{b=1}L_b\,g^{qk}\,v^s\,T^b_{qs}+\sum^n_{q=1}\sum^n_{r=1}
\sum^n_{s=1}\sum^n_{a=1}\sum^n_{j=1}L^r\,\times\\
\times\,g^{qk}\,\tilde\nabla_{\!q}g_{rs}\,T^s_{aj}\,v^a\,v^j
+\sum^n_{q=1}\sum^n_{s=1}\sum^n_{a=1}\sum^n_{j=1}
L_s\,g^{qk}\,\tilde\nabla_{\!q}T^s_{aj}\,v^a\,v^j
+\sum^n_{q=1}\sum^n_{j=1}\sum^n_{s=1}
2\,L_s\,\times\\
\times\,g^{qk}\,T^s_{qj}\,v^j-\sum^n_{q=1}\sum^n_{r=1}
\sum^n_{a=1}L^r\,g^{qk}\,T^a_{qr}\,L_a-\sum^n_{q=1}\sum^n_{r=1}
\sum^n_{a=1}\sum^n_{b=1}L_b\,g^{qk}\,T^b_{aq}\,v^a\,-\\
-\sum^n_{q=1}\sum^n_{m=1}\sum^n_{r=1}\sum^n_{s=1}
g^{mk}\,L_s\,\tilde\nabla_{\!m}T^s_{qr}\,L^r\,v^q.
\endgather
$$
Looking attentively at the above formula, we see that almost all
terms in right hand side do cancel each other. As a result we get
the following transformation rule:
$$
\hskip -2em
\alpha^k\to\alpha^k.
\tag9.17
$$
This means that first equation \thetag{5.14} is invariant under gauge
transformations \thetag{9.1}. In order to treat second equation
\thetag{5.14} we need to derive transformation rule for covector
fields $\boldsymbol\beta$ and $\boldsymbol\eta$ given by formulas
\thetag{7.31} and \thetag{7.32}. In \thetag{7.31} we have the entry
of second order covariant derivative $\nabla_{\!k}\!\nabla_{\!m}L_r$
and first order covariant derivative $\nabla_{\!k}F_r$. For 
$\nabla_{\!k}\!\nabla_{\!m}L_r$ we derive the following transformation
rule:
$$
\gathered
\nabla_{\!k}\!\nabla_{\!m}L_r\to\nabla_{\!k}\!\nabla_{\!m}L_r
-\sum^n_{a=1}\nabla_{\!k}T^a_{mr}\,L_a-\sum^n_{a=1}T^a_{mr}\,
\nabla_{\!k}L_a-\sum^n_{a=1}v^a\,\times\\
\times\,\sum^n_{b=1}\nabla_{\!k}g_{rb}\,T^b_{am}\,
-\sum^n_{a=1}\sum^n_{b=1}g_{rb}\,\nabla_{\!k}T^b_{am}\,v^a
-\sum^n_{c=1}T^c_{km}\,\nabla_{\!c}L_r+\sum^n_{c=1}T^c_{km}
\,\times\\
\times\,\sum^n_{a=1}T^a_{cr}\,L_a+\sum^n_{c=1}\sum^n_{a=1}
\sum^n_{b=1}T^c_{km}\,g_{rb}\,T^b_{ac}\,v^a-\sum^n_{c=1}
T^c_{kr}\nabla_{\!m}L_c+\sum^n_{c=1}T^c_{kr}\,\times\\
\times\,\sum^n_{a=1}T^a_{mc}\,L_a+\sum^n_{c=1}\sum^n_{a=1}
\sum^n_{b=1}T^c_{kr}\,g_{cb}\,T^b_{am}\,v^a-\sum^n_{c=1}
\sum^n_{e=1}T^c_{ek}\,\tilde\nabla_{\!c}\!\nabla_{\!m}L_r
\,\times\\
\times\,v^e+\sum^n_{c=1}\sum^n_{e=1}\sum^n_{a=1}v^e\,T^c_{ek}\,
\tilde\nabla_{\!c}T^a_{mr}\,L_a+\sum^n_{c=1}\sum^n_{e=1}
\sum^n_{a=1}v^e\,T^c_{ek}\,T^a_{mr}\,\tilde\nabla_{\!c}L_a\,+\\
+\sum^n_{a=1}\sum^n_{b=1}\sum^n_{c=1}\sum^n_{e=1}v^e\,T^c_{ek}
\,\tilde\nabla_{\!c}g_{rb}\,T^b_{am}\,v^a+\sum^n_{a=1}\sum^n_{b=1}
\sum^n_{c=1}\sum^n_{e=1}v^e\,T^c_{ek}\,g_{rb}\,\times\\
\times\,\tilde\nabla_{\!c}T^b_{am}\,v^a+\sum^n_{a=1}\sum^n_{b=1}
\sum^n_{e=1}v^e\,T^a_{ek}\,g_{rb}\,T^b_{am}.
\endgathered\quad
\tag9.18
$$
In a similar way for $\nabla_{\!k}F_r$ we derive transformation rule
which looks like:
$$
\gathered
\nabla_{\!k}F_r\to \nabla_{\!k}F_r+\sum^n_{s=1}\sum^n_{a=1}
\sum^n_{j=1}\nabla_{\!k}g_{rs}\,T^s_{aj}\,v^a\,v^j
+\sum^n_{s=1}\sum^n_{a=1}\sum^n_{j=1}g_{rs}\,\nabla_{\!k}T^s_{aj}
\,\times\\
\times\,v^a\,v^j-\sum^n_{c=1}T^c_{kr}
\,F_c-\sum^n_{c=1}T^c_{kr}\,\sum^n_{s=1}\sum^n_{a=1}
\sum^n_{j=1}g_{cs}\,T^s_{aj}\,v^a\,v^j-\sum^n_{c=1}\sum^n_{e=1}
v^e\,\times\\
\times\,T^c_{ke}\,\tilde\nabla_{\!c}F_r
-\sum^n_{c=1}\sum^n_{e=1}\sum^n_{s=1}\sum^n_{a=1}\sum^n_{j=1}
v^e\,T^c_{ke}\,\tilde\nabla_{\!c\,}g_{rs}\,T^s_{aj}\,v^a\,v^j
-\sum^n_{c=1}\sum^n_{e=1}v^e\,\times\\
\times\,\sum^n_{s=1}\sum^n_{a=1}\sum^n_{j=1}T^c_{ke}\,g_{rs}
\,\tilde\nabla_{\!c}T^s_{aj}\,v^a\,v^j-\sum^n_{e=1}\sum^n_{c=1}
\sum^n_{s=1}\sum^n_{a=1}2\,v^e\,T^c_{ke}\,g_{rs}
\,T^s_{ac}\,v^a.
\endgathered\quad
\tag9.19
$$
Now let's combine \thetag{9.1}, \thetag{9.2}, \thetag{9.4},
\thetag{9.5}, \thetag{9.6}, \thetag{9.7}, \thetag{9.8}, 
\thetag{9.9}, \thetag{9.10}, and two above formulas
\thetag{9.18} and \thetag{9.19}. Then for $\beta_k$ given
by formula \thetag{7.31} we obtain
$$
\allowdisplaybreaks
\gather
\beta_k\to\beta_k+\sum^n_{r=1}\sum^n_{c=1}\sum^n_{e=1}
v^r\,\nabla_{\!r}L^c\,T^e_{kc}\,L_e+\sum^n_{r=1}\sum^n_{c=1}
\sum^n_{e=1}v^r\,L^c\,\nabla_{\!r}T^e_{kc}\,L_e
+\sum^n_{r=1}v^r\,\times\\
\vspace{0pt plus 1pt}
\times\,\sum^n_{c=1}\sum^n_{e=1}L^c\,T^e_{kc}\,\nabla_{\!r}L_e
-\sum^n_{r=1}\sum^n_{m=1}v^r\,T^m_{rk}\,U_m-\sum^n_{r=1}
\sum^n_{m=1}\sum^n_{c=1}\sum^n_{e=1}v^r\,T^m_{rk}\,L^c\,T^e_{mc}
\,\times\\
\vspace{0pt plus 1pt}
\times\,L_e-\sum^n_{r=1}\sum^n_{m=1}\sum^n_{a=1}v^r\,v^a\,T^m_{ar}
\,\tilde\nabla_{\!m}U_k-\sum^n_{r=1}\sum^n_{m=1}\sum^n_{a=1}
\sum^n_{c=1}\sum^n_{e=1}v^r\,v^a\,T^m_{ar}\,\tilde\nabla_{\!m}L^c
\,L_e\,\times\\
\vspace{0pt plus 1pt}
\times\,T^e_{kc}
-\sum^n_{r=1}\sum^n_{m=1}\sum^n_{a=1}\sum^n_{c=1}\sum^n_{e=1}
v^r\,v^a\,T^m_{ar}\,L^c\,\tilde\nabla_{\!m}T^e_{kc}\,L_e
-\sum^n_{m=1}\sum^n_{c=1}\sum^n_{e=1}T^e_{kc}\,\tilde\nabla_{\!m}L_e
\,\times\\
\vspace{0pt plus 1pt}
\times\,\sum^n_{r=1}\sum^n_{a=1}v^r\,v^a\,T^m_{ar}\,L^c
{\vrule height 10pt depth 4pt}+\sum^n_{q=1}
\sum^n_{c=1}\sum^n_{r=1}F^q\,\tilde\nabla_{\!q}L^c\,T^r_{kc}\,L_r
+\sum^n_{q=1}\sum^n_{c=1}\sum^n_{r=1}F^q\,\tilde\nabla_{\!q}T^r_{kc}
\,\times\\
\vspace{0pt plus 1pt}
\times\,L^c\,L_r+\sum^n_{q=1}\sum^n_{c=1}\sum^n_{r=1}F^q\,L^c\,T^r_{kc}
\,\tilde\nabla_{\!q}L_r+\sum^n_{q=1}\sum^n_{i=1}\sum^n_{j=1}T^q_{ij}
\,v^i\,v^j\,\tilde\nabla_{\!q}U_k+\sum^n_{i=1}v^i\,\times\\
\vspace{0pt plus 1pt}
\times\,\sum^n_{c=1}\sum^n_{q=1}\sum^n_{r=1}\sum^n_{j=1}T^q_{ij}\,v^j\,
\tilde\nabla_{\!q}L^c\,T^r_{kc}\,L_r+\sum^n_{q=1}\sum^n_{c=1}\sum^n_{r=1}
\sum^n_{i=1}\sum^n_{j=1}T^q_{ij}\,v^i\,v^j\,L^c
\,\tilde\nabla_{\!q}T^r_{kc}\,\times\\
\vspace{0pt plus 1pt}
\times\,L_r+\sum^n_{q=1}\sum^n_{c=1}\sum^n_{r=1}\sum^n_{i=1}\sum^n_{j=1}
T^q_{ij}\,v^i\,v^j\,L^c\,T^r_{kc}\,\tilde\nabla_{\!q}L_r{\vrule height 10pt 
depth 4pt}+\sum^n_{q=1}
\sum^n_{r=1}\sum^n_{s=1}\sum^n_{c=1}T^c_{kq}\,L_c\,g^{rq}\,\times\\
\vspace{0pt plus 1pt}
\times\,v^s\,\nabla_{\!s}L_r+\sum^n_{q=1}\sum^n_{r=1}\sum^n_{s=1}
\sum^n_{c=1}\sum^n_{e=1}g_{qe}\,T^e_{ck}\,v^c\,g^{rq}\,v^s
\,\nabla_{\!s}L_r+\sum^n_{q=1}\sum^n_{r=1}\sum^n_{s=1}\nabla_{\!k}L_q
\,v^s\,\times\\
\vspace{0pt plus 1pt}
\times\,\sum^n_{a=1}g^{rq}\,T^a_{sr}\,L_a-\sum^n_{q=1}\sum^n_{r=1}
\sum^n_{s=1}\sum^n_{c=1}\sum^n_{a=1}T^c_{kq}\,L_c\,g^{rq}\,v^s\,
T^a_{sr}\,L_a-\sum^n_{q=1}\sum^n_{e=1}\sum^n_{r=1}g_{qe}\,\times\\
\vspace{0pt plus 1pt}
\times\,\sum^n_{s=1}\sum^n_{c=1}\sum^n_{a=1}T^e_{ck}\,v^c\,g^{rq}
\,v^s\,T^a_{sr}\,L_a+\sum^n_{q=1}\sum^n_{r=1}\sum^n_{s=1}\sum^n_{a=1}
\sum^n_{b=1}\nabla_{\!k}L_q\,g^{rq}\,v^s\,g_{rb}\,T^b_{as}\,\times\\
\vspace{0pt plus 1pt}
\times\,v^a-\sum^n_{q=1}\sum^n_{r=1}\sum^n_{s=1}\sum^n_{a=1}\sum^n_{b=1}
\sum^n_{c=1}T^c_{kq}\,L_c\,g^{rq}\,v^s\,g_{rb}\,T^b_{as}\,v^a-\sum^n_{q=1}
\sum^n_{r=1}\sum^n_{c=1}\sum^n_{e=1}g_{qe}\,\times\\
\vspace{0pt plus 1pt}
\times\,\sum^n_{s=1}\sum^n_{a=1}\sum^n_{b=1}T^e_{ck}\,v^c\,g^{rq}
\,v^s\,g_{rb}\,T^b_{as}\,v^a{\vrule height 10pt depth 4pt}+\sum^n_{q=1}
\sum^n_{r=1}\sum^n_{a=1}T^a_{kq}\,L_a\,g^{rq}\,F_r+\sum^n_{q=1}
\sum^n_{b=1}g_{qb}\,\times\\
\vspace{0pt plus 1pt}
\times\,\sum^n_{r=1}\sum^n_{a=1}T^b_{ak}\,v^a\,g^{rq}\,F_r-\sum^n_{q=1}
\sum^n_{r=1}\nabla_{\!k}L_q\,g^{rq}\,\sum^n_{s=1}\sum^n_{i=1}\sum^n_{j=1}
g_{rs}\,T^s_{ij}\,v^i\,v^j+\sum^n_{a=1}L_a\,\times\\
\vspace{0pt plus 1pt}
\times\,\sum^n_{q=1}\sum^n_{r=1}\sum^n_{s=1}\sum^n_{i=1}\sum^n_{j=1}
T^a_{kq}\,g^{rq}\,g_{rs}\,T^s_{ij}\,v^i\,v^j+\sum^n_{q=1}\sum^n_{r=1}
\sum^n_{a=1}\sum^n_{b=1}\sum^n_{s=1}\sum^n_{i=1}\sum^n_{j=1}g_{qb}
\,v^a\,\times\\
\vspace{0pt plus 1pt}
\times\,T^b_{ak}\,g^{rq}\,g_{rs}\,T^s_{ij}\,v^i\,v^j{\vrule height 10pt depth 
4pt}
+\sum^n_{r=1}\sum^n_{s=1}\sum^n_{a=1}\sum^n_{j=1}L^r\,\nabla_{\!k}g_{rs}
\,T^s_{aj}\,v^a\,v^j+\sum^n_{r=1}\sum^n_{s=1}\sum^n_{a=1}g_{rs}\,\times\\
\vspace{0pt plus 1pt}
\times\,\sum^n_{j=1}L^r\,\nabla_{\!k}T^s_{aj}\,v^a\,v^j-\sum^n_{r=1}
\sum^n_{c=1}L^r\,T^c_{kr}\,F_c-\sum^n_{r=1}\sum^n_{c=1}\sum^n_{s=1}
\sum^n_{a=1}\sum^n_{j=1}L^r\,T^c_{kr}\,g_{cs}\,v^a\,\times\\
\vspace{0pt plus 1pt}
\times\,T^s_{aj}\,v^j-\sum^n_{r=1}\sum^n_{c=1}\sum^n_{e=1}L^r\,v^e
\,T^c_{ke}\,\tilde\nabla_{\!c}F_r-\sum^n_{r=1}\sum^n_{c=1}
\sum^n_{e=1}\sum^n_{s=1}\sum^n_{a=1}\sum^n_{j=1}L^r\,v^e\,T^c_{ke}
\,T^s_{aj}\,\times\\
\vspace{0pt plus 1pt}
\times\,\tilde\nabla_{\!c\,}g_{rs}\,v^a\,v^j-\sum^n_{r=1}\sum^n_{c=1}
\sum^n_{e=1}\sum^n_{s=1}\sum^n_{a=1}\sum^n_{j=1}L^r\,v^e\,T^c_{ke}
\,g_{rs}\,\tilde\nabla_{\!c}T^s_{aj}\,v^a\,v^j-\sum^n_{r=1}\sum^n_{e=1}
L^r\,\times\\
\vspace{0pt plus 1pt}
\times\,\sum^n_{c=1}\sum^n_{s=1}\sum^n_{a=1}2\,v^e\,T^c_{ke}\,g_{rs}
\,T^s_{ac}\,v^a{\vrule height 10pt depth 4pt}-\sum^n_{m=1}\sum^n_{r=1}
\sum^n_{a=1}L^r\,v^m\,\nabla_{\!k}T^a_{mr}\,L_a-\sum^n_{m=1}\sum^n_{r=1}
L^r\,\times\\
\vspace{0pt plus 1pt}
\times\,\sum^n_{a=1}v^m\,T^a_{mr}\,\nabla_{\!k}L_a-\sum^n_{m=1}
\sum^n_{r=1}\sum^n_{a=1}\sum^n_{b=1}L^r\,v^m\,v^a\,\nabla_{\!k}g_{rb}
\,T^b_{am}-\sum^n_{m=1}\sum^n_{r=1}L^r\,v^m\,\times\\
\vspace{0pt plus 1pt}
\times\,\sum^n_{a=1}\sum^n_{b=1}g_{rb}\,\nabla_{\!k}T^b_{am}\,v^a
-\sum^n_{m=1}\sum^n_{r=1}\sum^n_{c=1}L^r\,v^m\,T^c_{km}\,\nabla_{\!c}L_r
+\sum^n_{m=1}\sum^n_{r=1}\sum^n_{c=1}\sum^n_{a=1}L^r\,\times\\
\vspace{0pt plus 1pt}
\times\,v^m\,T^c_{km}\,T^a_{cr}\,L_a+\sum^n_{m=1}\sum^n_{r=1}\sum^n_{c=1}
\sum^n_{a=1}\sum^n_{b=1}L^r\,v^m\,T^c_{km}\,g_{rb}\,T^b_{ac}\,v^a
-\sum^n_{m=1}\sum^n_{r=1}L^r\,v^m\,\times\\
\vspace{0pt plus 1pt}
\times\,\sum^n_{c=1}T^c_{kr}\nabla_{\!m}L_c+\sum^n_{m=1}\sum^n_{r=1}
\sum^n_{c=1}\sum^n_{a=1}L^r\,v^m\,T^c_{kr}\,T^a_{mc}\,L_a+\sum^n_{m=1}
\sum^n_{r=1}\sum^n_{c=1}L^r\,v^m\,T^c_{kr}\,\times\\
\vspace{0pt plus 1pt}
\times\,\sum^n_{a=1}\sum^n_{b=1}g_{cb}\,T^b_{am}\,v^a-\sum^n_{m=1}
\sum^n_{r=1}\sum^n_{c=1}\sum^n_{e=1}L^r\,v^m\,T^c_{ek}
\,\tilde\nabla_{\!c}\!\nabla_{\!m}L_r\,v^e+\sum^n_{m=1}\sum^n_{r=1}
L^r\,v^m\,\times\\
\vspace{0pt plus 1pt}
\times\,\sum^n_{c=1}\sum^n_{e=1}\sum^n_{a=1}v^e\,T^c_{ek}
\,\tilde\nabla_{\!c}T^a_{mr}\,L_a+\sum^n_{m=1}\sum^n_{r=1}\sum^n_{c=1}
\sum^n_{e=1}\sum^n_{a=1}L^r\,v^m\,v^e\,T^c_{ek}\,T^a_{mr}
\,\tilde\nabla_{\!c}L_a\,+\\
\vspace{0pt plus 1pt}
+\sum^n_{m=1}\sum^n_{r=1}\sum^n_{a=1}\sum^n_{b=1}\sum^n_{c=1}
\sum^n_{e=1}L^r\,v^m\,v^e\,T^c_{ek}\,\tilde\nabla_{\!c}g_{rb}
\,T^b_{am}\,v^a+\sum^n_{m=1}\sum^n_{r=1}\sum^n_{a=1}\sum^n_{b=1}
L^r\,v^m\,\times\\
\vspace{0pt plus 1pt}
\times\,\sum^n_{c=1}\sum^n_{e=1}v^e\,T^c_{ek}\,g_{rb}
\,\tilde\nabla_{\!c}T^b_{am}\,v^a+\sum^n_{m=1}\sum^n_{r=1}\sum^n_{a=1}
\sum^n_{b=1}\sum^n_{e=1}L^r\,v^m\,v^e\,T^a_{ek}\,g_{rb}\,T^b_{am}\,-
{\vrule height 10pt depth 4pt}\\
\vspace{0pt plus 1pt}
-\sum^n_{q=1}\sum^n_{r=1}\sum^n_{s=1}\sum^n_{c=1}\sum^n_{i=1}
\sum^n_{j=1}L^r\,\nabla_{\!k}L_q\,g^{sq}\,\tilde\nabla_{\!s}g_{rc}
\,T^c_{ij}\,v^i\,v^j-\sum^n_{q=1}\sum^n_{r=1}\sum^n_{s=1}
\sum^n_{c=1}\sum^n_{i=1}L^r\,\times\\
\vspace{0pt plus 1pt}
\times\,\sum^n_{j=1}\nabla_{\!k}L_q\,g^{sq}\,g_{rc}
\,\tilde\nabla_{\!s}T^c_{ij}\,v^i\,v^j-\sum^n_{q=1}\sum^n_{r=1}
\sum^n_{s=1}\sum^n_{c=1}\sum^n_{j=1}2\,L^r\,\nabla_{\!k}L_q\,g^{sq}
\,g_{rc}\,T^c_{sj}\,v^j\,+\\
\vspace{0pt plus 1pt}
+\sum^n_{q=1}\sum^n_{r=1}\sum^n_{s=1}\sum^n_{a=1}L^r\,T^a_{kq}\,L_a
\,g^{sq}\,\tilde\nabla_{\!s}F_r+\sum^n_{q=1}\sum^n_{r=1}\sum^n_{s=1}
\sum^n_{a=1}\sum^n_{c=1}\sum^n_{i=1}\sum^n_{j=1}L^r\,T^a_{kq}\,L_a
\,\times\\
\vspace{0pt plus 1pt}
\times\,g^{sq}\,\tilde\nabla_{\!s}g_{rc}\,T^c_{ij}\,v^i\,v^j+\sum^n_{q=1}
\sum^n_{r=1}\sum^n_{s=1}\sum^n_{a=1}\sum^n_{c=1}\sum^n_{i=1}\sum^n_{j=1}
L^r\,T^a_{kq}\,L_a\,g^{sq}\,g_{rc}\,\tilde\nabla_{\!s}T^c_{ij}\,v^i\,v^j
\,+\\
\vspace{0pt plus 1pt}
+\sum^n_{q=1}\sum^n_{r=1}\sum^n_{s=1}\sum^n_{a=1}\sum^n_{c=1}
\sum^n_{j=1}2\,L^r\,T^a_{kq}\,L_a\,g^{sq}\,g_{rc}\,T^c_{sj}\,v^j
+\sum^n_{q=1}\sum^n_{r=1}\sum^n_{s=1}\sum^n_{a=1}\sum^n_{b=1}
L^r\,g_{qb}\,\times\\
\vspace{0pt plus 1pt}
\times\,T^b_{ak}\,v^a\,g^{sq}\,\tilde\nabla_{\!s}F_r+\sum^n_{q=1}
\sum^n_{r=1}\sum^n_{s=1}\sum^n_{a=1}\sum^n_{b=1}\sum^n_{c=1}
\sum^n_{i=1}\sum^n_{j=1}L^r\,g_{qb}\,T^b_{ak}\,v^a\,g^{sq}\,
\tilde\nabla_{\!s}g_{rc}\,T^c_{ij}\,\times\\
\vspace{0pt plus 1pt}
\times\,v^i\,v^j+\sum^n_{q=1}\sum^n_{r=1}\sum^n_{s=1}\sum^n_{a=1}
\sum^n_{b=1}\sum^n_{c=1}\sum^n_{i=1}\sum^n_{j=1}L^r\,g_{qb}\,T^b_{ak}
\,v^a\,g^{sq}\,g_{rc}\,\tilde\nabla_{\!s}T^c_{ij}\,v^i\,v^j
+\sum^n_{r=1}2\,\times\\
\vspace{0pt plus 1pt}
\times\,\sum^n_{q=1}\sum^n_{s=1}\sum^n_{a=1}\sum^n_{b=1}\sum^n_{c=1}
\sum^n_{j=1}L^r\,g_{qb}\,T^b_{ak}\,v^a\,g^{sq}\,g_{rc}\,T^c_{sj}\,v^j
{\vrule height 10pt depth 4pt}+\sum^n_{q=1}\sum^n_{m=1}\sum^n_{r=1}
\sum^n_{s=1}\sum^n_{c=1}L^r\,\times\\
\vspace{0pt plus 1pt}
\times\,\nabla_{\!k}L_q\,g^{sq}\,v^m\,\tilde\nabla_{\!s}T^c_{mr}\,L_c
+\sum^n_{q=1}\sum^n_{m=1}\sum^n_{r=1}\sum^n_{s=1}\sum^n_{c=1}
L^r\,\nabla_{\!k}L_q\,g^{sq}\,v^m\,T^c_{mr}\,g_{cs}+\sum^n_{r=1}
L^r\,\times\\
\vspace{0pt plus 1pt}
\times\,\sum^n_{q=1}\sum^n_{m=1}\sum^n_{s=1}\sum^n_{c=1}\sum^n_{e=1}
\nabla_{\!k}L_q\,g^{sq}\,v^m\,\tilde\nabla_{\!s}g_{re}\,T^e_{cm}\,v^c
+\sum^n_{q=1}\sum^n_{m=1}\sum^n_{r=1}\sum^n_{s=1}\sum^n_{c=1}
\sum^n_{e=1}L^r\,\times\\
\vspace{0pt plus 1pt}
\times\,\nabla_{\!k}L_q\,g^{sq}\,v^m\,g_{re}\,\tilde\nabla_{\!s}T^e_{cm}
\,v^c+\sum^n_{q=1}\sum^n_{m=1}\sum^n_{r=1}\sum^n_{s=1}\sum^n_{e=1}L^r
\,\nabla_{\!k}L_q\,g^{sq}\,v^m\,g_{re}\,T^e_{sm}\,+\\
\vspace{0pt plus 1pt}
+\sum^n_{q=1}\sum^n_{m=1}\sum^n_{r=1}\sum^n_{s=1}\sum^n_{a=1}L^r
\,T^a_{kq}\,L_a\,g^{sq}\,v^m\,\tilde\nabla_{\!s}\!\nabla_{\!m}L_r
-\sum^n_{q=1}\sum^n_{m=1}\sum^n_{r=1}\sum^n_{s=1}\sum^n_{a=1}
\sum^n_{c=1}L^r\,\times\\
\vspace{0pt plus 1pt}
\times\,T^a_{kq}\,L_a\,g^{sq}\,v^m\,\tilde\nabla_{\!s}T^c_{mr}
\,L_c-\sum^n_{q=1}\sum^n_{m=1}\sum^n_{r=1}\sum^n_{s=1}\sum^n_{a=1}
\sum^n_{c=1}L^r\,T^a_{kq}\,L_a\,g^{sq}\,v^m\,T^c_{mr}\,g_{cs}\,-\\
\vspace{0pt plus 1pt}
-\sum^n_{q=1}\sum^n_{m=1}\sum^n_{r=1}\sum^n_{s=1}\sum^n_{a=1}
\sum^n_{c=1}\sum^n_{e=1}L^r\,T^a_{kq}\,L_a\,g^{sq}\,v^m
\,\tilde\nabla_{\!s}g_{re}\,T^e_{cm}\,v^c-\sum^n_{q=1}\sum^n_{m=1}
\sum^n_{r=1}L^r\,\times\\
\vspace{0pt plus 1pt}
\times\,\sum^n_{s=1}\sum^n_{a=1}\sum^n_{c=1}\sum^n_{e=1}T^a_{kq}
\,L_a\,g^{sq}\,v^m\,g_{re}\,\tilde\nabla_{\!s}T^e_{cm}\,v^c
-\sum^n_{q=1}\sum^n_{m=1}\sum^n_{r=1}\sum^n_{s=1}\sum^n_{a=1}
\sum^n_{e=1}L^r\,T^a_{kq}\,\times\\
\vspace{0pt plus 1pt}
\times\,L_a\,g^{sq}\,v^m\,g_{re}\,T^e_{sm}+\sum^n_{q=1}\sum^n_{m=1}
\sum^n_{r=1}\sum^n_{s=1}\sum^n_{a=1}\sum^n_{b=1}L^r\,g_{qb}\,T^b_{ak}
\,v^a\,g^{sq}\,v^m\,\tilde\nabla_{\!s}\!\nabla_{\!m}L_r\,-\\
\vspace{0pt plus 1pt}
-\sum^n_{q=1}\sum^n_{m=1}\sum^n_{r=1}\sum^n_{s=1}\sum^n_{a=1}
\sum^n_{b=1}\sum^n_{c=1}L^r\,g_{qb}\,T^b_{ak}\,v^a\,g^{sq}\,v^m
\,\tilde\nabla_{\!s}T^c_{mr}\,L_c-\sum^n_{q=1}\sum^n_{m=1}
\sum^n_{r=1}L^r\,\times\\
\vspace{0pt plus 1pt}
\times\,\sum^n_{s=1}\sum^n_{a=1}\sum^n_{b=1}\sum^n_{c=1}g_{qb}
\,T^b_{ak}\,v^a\,g^{sq}\,v^m\,T^c_{mr}\,g_{cs}-\sum^n_{q=1}
\sum^n_{m=1}\sum^n_{r=1}\sum^n_{s=1}\sum^n_{a=1}\sum^n_{b=1}
L^r\,g_{qb}\,T^b_{ak}\,\times\\
\vspace{0pt plus 1pt}
\times\,\sum^n_{c=1}\sum^n_{e=1}v^a\,g^{sq}\,v^m
\,\tilde\nabla_{\!s}g_{re}\,T^e_{cm}\,v^c-\sum^n_{q=1}\sum^n_{m=1}
\sum^n_{r=1}\sum^n_{s=1}\sum^n_{a=1}\sum^n_{b=1}L^r\,g_{qb}\,T^b_{ak}
\,v^a\,g^{sq}\,\times\\
\vspace{0pt plus 1pt}
\times\,v^m\,\sum^n_{c=1}\sum^n_{e=1}g_{re}
\,\tilde\nabla_{\!s}T^e_{cm}\,v^c-\sum^n_{q=1}\sum^n_{m=1}\sum^n_{r=1}
\sum^n_{s=1}\sum^n_{a=1}\sum^n_{b=1}\sum^n_{e=1}L^r\,g_{qb}\,T^b_{ak}
\,v^a\,g^{sq}\,v^m\,\times\\
\vspace{0pt plus 1pt}
\times\,g_{re}\,T^e_{sm}{\vrule height 10pt depth 4pt}
-\sum^n_{r=1}\sum^n_{s=1}\sum^n_{m=1}\nabla_{\!m}T^s_{kr}\,v^m
\,L^r\,L_s+\sum^n_{r=1}\sum^n_{s=1}\sum^n_{m=1}\nabla_{\!k}T^s_{mr}
\,v^m\,L^r\,L_s\,+\\
\vspace{0pt plus 1pt}
+\sum^n_{r=1}\sum^n_{s=1}\sum^n_{m=1}\sum^n_{b=1}\sum^n_{a=1}v^a
\,T^b_{ka}\,D^s_{mrb}\,v^m\,L^r\,L_s-\sum^n_{r=1}\sum^n_{s=1}
\sum^n_{m=1}\sum^n_{b=1}\sum^n_{a=1}v^a\,T^b_{ma}\,v^m\,\times\\
\vspace{0pt plus 1pt}
\times\,D^s_{krb}\,L^r\,L_s-\sum^n_{r=1}\sum^n_{s=1}\sum^n_{m=1}
\sum^n_{a=1}T^s_{ma}\,T^a_{kr}\,v^m\,L^r\,L_s+\sum^n_{r=1}\sum^n_{s=1}
\sum^n_{m=1}\sum^n_{a=1}T^s_{ka}\,T^a_{mr}\,\times\\
\vspace{0pt plus 1pt}
\times\,v^m\,L^r\,L_s-\sum^n_{r=1}\sum^n_{s=1}\sum^n_{m=1}\sum^n_{b=1}
\sum^n_{a=1}v^a\,T^b_{ka}\,\tilde\nabla_{\!b}T^s_{mr}\,v^m\,L^r\,L_s
+\sum^n_{r=1}\sum^n_{s=1}\sum^n_{m=1}\sum^n_{a=1}v^a\,\times\\
\vspace{0pt plus 1pt}
\times\,\sum^n_{b=1}T^b_{ma}\,\tilde\nabla_{\!b}T^s_{kr}\,v^m
\,L^r\,L_s{\vrule height 10pt depth 4pt}
-\sum^n_{q=1}\sum^n_{r=1}\sum^n_{s=1}\sum^n_{m=1}\sum^n_{c=1}
\sum^n_{a=1}g^{cq}\,T^a_{kq}\,L_a\,D^s_{mrc}\,v^m\,L^r\,L_s\,-\\
\vspace{0pt plus 1pt}
-\sum^n_{q=1}\sum^n_{r=1}\sum^n_{s=1}\sum^n_{m=1}\sum^n_{a=1}
\sum^n_{c=1}\sum^n_{b=1}g^{cq}\,g_{qb}\,T^b_{ak}\,v^a\,D^s_{mrc}
\,v^m\,L^r\,L_s-\sum^n_{q=1}\sum^n_{r=1}\sum^n_{c=1}g^{cq}\,\times\\
\vspace{0pt plus 1pt}
\times\,\sum^n_{m=1}\sum^n_{s=1}\nabla_{\!k}L_q
\,\tilde\nabla_{\!c}T^s_{rm}\,v^m\,L^r\,L_s+\sum^n_{q=1}\sum^n_{r=1}
\sum^n_{s=1}\sum^n_{m=1}\sum^n_{c=1}\sum^n_{a=1}g^{cq}\,T^a_{kq}
\,L_a\,\tilde\nabla_{\!c}T^s_{rm}\,\times\\
\vspace{0pt plus 1pt}
\times\,v^m\,L^r\,L_s+\sum^n_{q=1}\sum^n_{r=1}\sum^n_{s=1}\sum^n_{m=1}
\sum^n_{c=1}\sum^n_{a=1}\sum^n_{b=1}g^{cq}\,g_{qb}\,T^b_{ak}\,v^a
\,\tilde\nabla_{\!c}T^s_{rm}\,v^m\,L^r\,L_s
{\vrule height 10pt depth 4pt}
+\sum^n_{i=1}v^i\,\times\\
\vspace{0pt plus 1pt}
\times\,\sum^n_{r=1}\sum^n_{s=1}\sum^n_{m=1}\sum^n_{a=1}\sum^n_{b=1}
\sum^n_{j=1}g^{am}\,D^s_{rka}\,g_{mb}\,T^b_{ij}\,v^j\,L^r\,L_s
-\sum^n_{r=1}\sum^n_{s=1}\sum^n_{m=1}\sum^n_{a=1}g^{am}\,\times\\
\vspace{0pt plus 1pt}
\times\,\tilde\nabla_{\!a}T^s_{kr}\,F_m\,L^r\,L_s-\sum^n_{r=1}
\sum^n_{s=1}\sum^n_{m=1}\sum^n_{a=1}\sum^n_{b=1}\sum^n_{i=1}
\sum^n_{j=1}g^{am}\,\tilde\nabla_{\!a}T^s_{kr}\,g_{mb}\,T^b_{ij}
\,v^i\,v^j\,L^r\,L_s.\\
\endgather
$$
Totally in right hand side of the above formula we have 97 terms
apart from $\beta_k$. Lets denote them by $T[i]$. Then formula
can be written as
$$
\hskip -2em
\beta_k\to\beta_k+\sum^{97}_{i=1}T[i].
\tag9.20
$$
Terms $T[1],\,\ldots,\,T[97]$ are divided into 11 groups separated
from each other by vertical bars. These groups represent contributions
from 11 terms in formula \thetag{7.31}:
$$
97=9+7+8+5+8+16+11+17+8+5+3.
$$
One can find that most of terms $T[1],\,\ldots,\,T[97]$ in formula
\thetag{9.20} do cancel each other. In the following table we
collect mutually canceling terms:
$$
\xalignat 3
&T[2]+T[82]=0,&&T[3]+T[45]=0,&&T[5]+T[43]=0,\\
&T[6]+T[13]=0,&&T[7]+T[14]=0,&&T[8]+T[15]=0,\\
&T[9]+T[16]=0,&&T[11]+T[96]=0,&&T[12]+T[32]=0,\\
&T[21]+T[44]=0,&&T[22]+T[27]=0,&&T[23]+T[28]=0,\\
&T[24]+T[29]=0,&&T[30]+T[40]=0,&&T[31]+T[41]=0,\\
&T[33]+T[47]=0,&&T[34]+T[61]=0,&&T[35]+T[51]=0,\\
&T[36]+T[52]=0,&&T[37]+T[64]=0,&&T[38]+T[83]=0,\\
&T[39]+T[66]=0,&&T[46]+T[86]=0,&&T[48]+T[76]=0,\\
&T[49]+T[77]=0,&&T[50]+T[78]=0,&&T[53]+T[81]=0,\\
&T[54]+T[67]=0,&&T[55]+T[68]=0,&&T[58]+T[73]=0,\\
&T[59]+T[74]=0,&&T[62]+T[79]=0,&&T[63]+T[80]=0,\\
&T[65]+T[92]=0,&&T[71]+T[93]=0,&&T[72]+T[87]=0,\\
&T[84]+T[91]=0,&&T[85]+T[95]=0,&&T[88]+T[94]=0,\\
\vspace{2ex}
&\qquad T[19]+T[56]+T[69]=0,\kern -10em
&&\qquad\qquad\quad T[20]+T[60]+T[75]=0,\kern -10em\\
&\qquad T[89]+T[97]=0,
&&\qquad\qquad\quad T[4]+T[18]+T[26]+T[42]=0.\kern -10em
\endxalignat
$$
After all above cancellations only 7 terms $T[i]$ survive in
formula \thetag{9.20} for $\beta_k$:
$$
\beta_k\to\beta_k+T[1]+T[10]+T[17]+T[25]+T[57]+T[70]+T[90].
\quad
\tag9.21
$$
Now we can write \thetag{9.21} in explicit form. This form is rather
observable one:
$$
\gathered
\beta_k\to\beta_k+\sum^n_{r=1}\sum^n_{c=1}\sum^n_{e=1}v^r
\,\nabla_{\!r}L^c\,T^e_{kc}\,L_e+\sum^n_{q=1}\sum^n_{c=1}
\sum^n_{r=1}F^q\,\tilde\nabla_{\!q}L^c\,T^r_{kc}\,L_r\,+\\
+\sum^n_{q=1}\sum^n_{r=1}\sum^n_{s=1}\sum^n_{c=1}T^c_{kq}
\,L_c\,g^{rq}\,v^s\,\nabla_{\!s}L_r+\sum^n_{q=1}\sum^n_{r=1}
\sum^n_{a=1}T^a_{kq}\,L_a\,g^{rq}\,F_r\,+\\
+\sum^n_{q=1}\sum^n_{r=1}\sum^n_{s=1}\sum^n_{a=1}L^r\,T^a_{kq}
\,L_a\,g^{sq}\,\tilde\nabla_{\!s}F_r+\sum^n_{q=1}\sum^n_{r=1}
\sum^n_{m=1}\sum^n_{s=1}\sum^n_{a=1}L^r\,T^a_{kq}\,L_a\,\times\\
\times\,g^{sq}\,v^m\,\tilde\nabla_{\!s}\!\nabla_{\!m}L_r
-\sum^n_{q=1}\sum^n_{r=1}\sum^n_{s=1}\sum^n_{m=1}\sum^n_{c=1}
\sum^n_{a=1}g^{cq}\,T^a_{kq}\,L_a\,D^s_{mrc}\,v^m\,L^r\,L_s.\\
\endgathered
\tag9.22
$$
Looking attentively at the above formula \thetag{9.22}, one can
find that it can be further transformed. Namely, one can extract
common factors in all terms:
$$
\gathered
\beta_k\to\beta_k+\sum^n_{e=1}\sum^n_{q=1}T^e_{kq}\,L_e\left(
\,\shave{\sum^n_{r=1}}v^r\,\,\nabla_{\!r}L^q+\shave{\sum^n_{r=1}}
F^r\,\tilde\nabla_{\!r}L^q\,+\right.\\
+\shave{\sum^n_{s=1}}\shave{\sum^n_{r=1}}\shave{\sum^n_{c=1}}
L^r\,g^{cq}\,v^s\,\tilde\nabla_{\!c}\!\nabla_{\!s}L_r
+\sum^n_{r=1}
g^{rq}\,F_r+\sum^n_{r=1}\sum^n_{s=1}L^r\,g^{sq}\,
\tilde\nabla_{\!s}F_r\,+\\
\left.+\shave{\sum^n_{r=1}}\shave{\sum^n_{s=1}}g^{rq}\,v^s\,
\nabla_{\!s}L_r-\shave{\sum^n_{m=1}}\shave{\sum^n_{r=1}}
\shave{\sum^n_{s=1}}\shave{\sum^n_{c=1}}g^{mq}\,L_s\,D^s_{crm}
\,L^r\,v^c\!\right)\!.
\endgathered
\tag9.23
$$
Due to the equality \thetag{7.8} we can derive the following relationship
for $U_r$ and $F_r$:
$$
\hskip -2em
\sum^n_{r=1}g^{rq}\left(\!F_r+\shave{\sum^n_{s=1}}v^s\,\nabla_sL_r
\!\right)=\sum^n_{r=1}g^{rq}\left(\!U_r+\shave{\sum^n_{s=1}}L^s\,
\nabla_rL_s\!\right)\!.
\tag9.24
$$
Combining \thetag{9.24} with \thetag{9.23} and applying formula
\thetag{7.19}, we can provide substantial simplification of formula
\thetag{7.23}. Now it looks like
$$
\hskip -2em
\beta_k\to\beta_k+\sum^n_{e=1}\sum^n_{q=1}T^e_{kq}\,L_e\,\alpha^q.
\tag9.25
$$
Let's remember formula \thetag{7.32} for $\eta_k$. Applying
\thetag{9.6}, \thetag{9.17}, and \thetag{9.25} to \thetag{7.32}
and using formula \thetag{6.4} for projector components, we derive
$$
\hskip -2em
\eta_k\to\eta_k+\sum^n_{e=1}\sum^n_{q=1}\sum^n_{s=1}T^e_{kq}\,L_e
\,P^q_s\,\alpha^s.
\tag9.26
$$
If we recall first equation \thetag{5.14}, we see that \thetag{9.26}
is equivalent to
$$
\pagebreak
\hskip -2em
\eta_k\to\eta_k
\tag9.27
$$
This means that covector field $\boldsymbol\eta$ is invariant
under gauge transformations \thetag{9.1}. Due to \thetag{9.17}
and \thetag{9.27} both weak normality equations \thetag{5.14}
are invariant under gauge transformations \thetag{9.1}. This
result can be stated as a theorem.
\proclaim{Theorem 9.2} Weak normality equations \thetag{5.14}
transformed to $\bold v$-representation are invariant under
gauge transformations \thetag{9.1}.
\endproclaim
Theorems~9.1 and 9.2 mean that normality equations \thetag{5.14},
\thetag{5.18}, \thetag{5.19}, and \thetag{5.20}, when applied to
Newtonian dynamical system \thetag{1.1}, do not actually depend
on connection components $\Gamma^k_{ij}$. Therefore we can
substitute 
$$
\xalignat 4
&\Gamma^k_{ij}=0,&&F^i=\Phi^i,&&\nabla_i=\frac{\partial}
{\partial x^i},&&\tilde\nabla_i=\frac{\partial}{\partial v^i}
\endxalignat
$$
and write normality equations in connection free form. In such
form each term of these equations has no separate tensorial
interpretation and one should find invariant (coordinate-free)
interpretation for the equations in whole. However, this is
separate problem which will be studied later.
\head
10. Acknowledgements.
\endhead
     This work is supported by grant from Russian Fund for Basic
Research (project 00-01-00068, coordinator of project
Ya\.~T.~Sultanaev), and by grant
from Academy of Sciences of the Republic Bashkortostan (coordinator
N.~M.~Asadullin). I am grateful to these organizations for financial
support.

\enddocument
\end